%% file: stratified.tex
\input abbrev

\ref\AGJ{J. Arms, M. Gotay and G. Jennings, Geometric and algebraic reduction for singular momentum maps, Adv. in Math. 79 (1990) 43-103.}

\ref\Atiyah{M. Atiyah, Convexity and commuting Hamiltonians, Bull. London Math. Soc. 14 (1982) 1-15.}



\ref\Delzant{T. Delzant, Hamiltoniens p\'eriodiques et image convexe de l'application moment, Bull. Soc. Math.
France 116 (1988) 315-339.}

\ref\GM{M. Goresky and R. MacPherson, Intersection homology theory,
Topology 19 (1980) 135-162.}

\ref\Guillemin{V. Guillemin, Moment maps and Combinatorial Invariants of Hamiltonian
$T^n$-spaces, Progress in Mathematics vol. 122, Birkhauser, Boston 1994.}

\ref\GS{V. Guillemin and S. Sternberg, Convexity properties of the moment mapping,
Invent. Math. 67 (1982) 491-513.}

\ref\GSI{V. Guillemin and S. Sternberg, A normal form for the moment map, in: Differential Geometric Methods in Mathematical Physics, S. Sternberg, editor, Reidel Publishing Company, Dordrecht, 1984.}

\ref\Kirwan{F. Kirwan, Convexity properties of the moment map. III, Invent. Math. 77 (1984) 547-552.}

\ref\LermanTolman{E. Lerman and S. Tolman, Hamiltonian torus actions on symplectic orbifolds and toric varieties. Trans. Amer. Math. Soc. 349 (1997), no. 10, 4201--4230.} 

\ref\Ma{C.-M. Marle, Mod\`ele d'action hamiltonienne d'un groupe de Lie sur une vari\'et\'e symplectique, in: Rendiconti del Seminario Matematico, Universit\`a e Politecnico, Torino 43 (1985) 227-251.}

\ref\MW{J. Marsden and A. Weinstein,
Reduction of symplectic manifolds with symmetry,
Rep. Math. Phys. Vol. 5 No. 1 (1974) 121-130.}

\ref\OrtegaRatiu{J.-P. Ortega and T.S. Ratiu, Momentum maps and Hamiltonian reduction. Prog. in Math., 222. BirkhŠuser Boston, Inc., Boston, MA, 2004.} 

\ref\SL{R. Sjamaar and E. Lerman,
Stratified symplectic spaces and reduction, Annals Math. 134 (1991) 375-422.}

\ref\Sternberg{S. Sternberg, On minimal coupling and the symplectic mechanics of a classical particle in the presence of a Yang-Mills field, Proc. Nat. Acad. Sci. USA 74 (1977) 5253-5254.}

\ref\Weinstein{A. Weinstein, A universal phase space for particles in Yang-Mills fields, Lett. Math. Phys. 2 (1978) 417-420.}

\nopagenumbers

\centerline{\titlefont A Convexity Theorem and}
\centerline{\titlefont Reduced Delzant Spaces}

\bs
\centerline{Bong H. Lian and Bailin Song}
\bs

ABSTRACT.  The convexity theorem of Atiyah and Guillemin-Sternberg says that any connected compact manifold with Hamiltonian torus action has a moment map whose image is the convex hull of the image of the fixed point set. Sjamaar-Lerman proved that the Marsden-Weinstein reduction of a connected Hamitonian $G$-manifold is a stratified symplectic space. Suppose $1\ra A\ra G\ra T\ra 1$ is an exact sequence of compact Lie groups and $T$ is a torus. Then the reduction of a Hamiltonian $G$-manifold with respect to $A$ yields a Hamiltonian $T$-space. 
We show that if the $A$-moment map is proper, then the convexity theorem holds for such a Hamiltonian $T$-space, even when it is singular. We also prove that if, furthermore, the $T$-space has dimension $2dim~T$ and $T$ acts effectively, then the moment polytope is sufficient to essentially distinguish their homeomorphism type, though not their diffeomorphism types. This generalizes a theorem of Delzant in the smooth case.

\bs

\np
\headline{\ifodd\pageno\rightheadline\else\leftheadline\fi}
\def\rightheadline{\tenrm\hfil Convexity Theorem and Reduced Delzant Spaces 
\hfil\folio}
\def\leftheadline{\tenrm\folio\hfil Bong H. Lian \& Bailin Song \hfil}

\newsec{Introduction}

Let $G$ be a compact Lie group, and $(M,\omega)$ a connected symplectic manifold with a Hamiltonian $G$-action and a moment map $J$. Then
the Marsden-Weinstein reduced space $M_0=J^{-1}(0)/G$ can in general be a singular space \MW.
Arms-Gotay-Jennings \AGJ~ introduced the smooth structure $C^\infty(M_0)=C^\infty(M)^G/I^G$ on $M_0$, where $I^G$ is the ideal of invariant functions vanishing on $J^{-1}(0)$,  and showed that
$C^\infty(M_0)$ inherits a Poisson structure from $C^\infty(M)$. 
Sjamaar-Lerman \SL~ proved that $M_0$ has the structure of a stratified space in the sense of Goresky-MacPherson \GM. Moreover, the stratified space $M_0$ is symplectic in the sense that the Poisson structure is compatible with the symplectic structure on each stratum of $M_0$ (Definition 1.12 \SL).
Another remarkable result in \SL~ says that  if $J$ is proper, the reduced space $M_0$ always contains a unique open connected dense stratum. Clearly a stratified symplectic space, in general, need not have this property. Thus Sjamaar-Lerman's result can be thought of as a kind of obstruction: if a stratified symplectic space has no connected open dense stratum then it cannot be the reduction of a smooth Hamiltonian manifold.

Throughout this paper $T$ denotes a torus.

\definition{A stratified space $X$ is called a reduced $T$-space if there is an exact sequence
of compact Lie groups  $1\ra A{\br i\over\ra} G\ra T\ra 1$ and a connected Hamiltonian $G$-manifold $M$ with a moment map $J$, such that $X=(i^*\circ J)^{-1}(0)/A$, where the $A$-moment map $i^*\circ J:M\ra\ga^*$ is assumed proper.}

By Theorem 4.4 \SL,  a reduced $T$-space is a stratified symplectic $T$-space with a moment map
$J_T$ induced by $J$. In section 2, we prove$^1$\footnote{}{After the completion of this paper,
in a joint work with Dong Wang we have extended this theorem to the case when $T$ is nonabelian, hence generalizing a theorem of Kirwan \Kirwan~to reduced stratified $T$ spaces.}

\theorem{Let $X$ be a reduced $T$-space. Then the induced moment map $J_T$ has the following properties:
\item{i.} the level sets of $J_T$ are connected;
\item{ii.} the image of $J_T$ is convex;
\item{iii.} the image of $J_T$ is the convex hull of the image of the
fixed points of the $T$ action.}
\thmlab\TheoremI

Part i. generalizes a theorem of Atiyah, while ii.-iii. generalizes
the convexity theorem of Atiyah \Atiyah~ and Guillemin-Sternberg \GS.
The convexity of $J_T(X)$ provides a new obstruction: if a stratified symplectic $T$-space (in the sense of \SL) has non-convex moment image $J_T(X)$ then it cannot be the reduction of a Hamiltonian manifold. There are examples that show that this obstruction is independent of the Sjamaar-Lerman obstruction.

\definition{A reduced $T$-space $X$ is called a reduced Delzant $T$-space if
$dim~X=2dim~T$ and $T$ acts effectively on the dense open stratum of $X$.}

\definition{We say that an $n$ dimensional polytope in $\gt^*=Lie(T^n)^*$ is rational (or $T$-rational) if every facet has a normal vector $u\in\gt$ which is rational, i.e. it generates a closed subgroup in $T$.}
\thmlab\DefRationality

\theorem{
The moment polytope of a reduced Delzant $T$-space is a rational polytope. 
Every rational polytope can be realized as the moment polytope of a reduced Delzant $T$-space which has the structure of a complete toric variety.}
\thmlab\TheoremII

\theorem{Let $X$ be a reduced Delzant $T$-space. Assume that the stabilizer of each point $x$ in $X$ is connected.  Then $J_T(X)$ determines the homeomorphism type of $X$. In fact if $(X,J_T)$,
$(X',J'_T)$ are two reduced Delzant $T$-spaces with the same moment polytope in $\gt^*$, then there exists a $T$-equivariant homeomorphism $\varphi:X\ra X'$ such that $J_T=J'_T\circ\varphi$.}
\thmlab\TheoremIII

These two theorems generalize results of Delzant \Delzant~ in the smooth case and of Lerman-Tolman \LermanTolman~the orbifold case;  see Guillemin's book \Guillemin~ for a more extensive review. 
The second theorem implies that the homeomorphism types of those reduced Delzant $T$-spaces are classified by their moment polytopes, and that all of them are realized by toric varieties. Note however that one cannot hope to recover
the smooth structure of a Delzant $T$-space in this generality. For example, consider the projective line $\P^1$ and $\P^1/\Z_n$, where $\Z_n$ is the cyclic subgroup of order $n$ in $S^1$
acting by standard rotations on $\P^1$. Both can be realized as reduced Delzant $S^1$-spaces.
By making a suitable choice of symplectic structures, we can also make their moment polytopes  equal.
Yet the two spaces can have different smooth structures. The condition on the stabilizers in $X$
is a technical assumption which we conjecture to be superfluous. When $X$ is smooth,
this is a consequence of the equivariant Darboux theorem. We will prove the Theorem \TheoremII~in sections 3-4, and Theorem \TheoremIII~ in section 5.

There is a construction in algebraic geometry that realizes toric varieties as a kind of ``categorical quotients''. 
It should be emphasized that Theorem \TheoremII~above does not rely on results in algebraic geometry, and is not a consequence of the categorical quotient construction. What the theorem shows is that those toric varieties can also be realized {\it topologically} as a symplectic stratified space a la Delzant. There has been assertions in the literature that mistakenly claim that only a {\it projective} toric variety can have the structure of a Delzant space. One should not confuse a projective structure of a variety, which is a feature of the variety's  complex structure, with its symplectic structure as a stratifed space. In fact, as an immediate consequence of Delzant's result, any smooth complete toric variety, projective or not, is {\it diffeomorphic} (in fact symplectically if a suitable Kahler structure is chosen) to a Delzant space. Delzant's diffeomorphism is about the symplectic diffeomorphism classes, and not about the complex algebraic isomorphism classes. Likewise, our result is not about the latter. But rather it is about the homeomorphism classes of symplectic stratified spaces. The reader should be cautioned to avoid this confusion here.

We now comment on the main ideas of our proofs. 
The proof of Theorem \TheoremI~follows a strategy that is parallel to that in the smooth case. We want to show that the intersection of image of the moment map with any rational line is connected. At a crucial point, we need to use the so-called local normal form of the moment map, discovered by Marle \Ma~ and Guillemin-Sternberg \GSI, and in the form we will use, developed by Sjamaar-Lerman \SL. This is needed to establish that the fixed point set is a stratified space having only finitely many components.
The image of the moment map is then recovered as the the convex hull of the image of fixed point set, as in the smooth case. Theorems \TheoremII~and \TheoremIII~are a bit more delicate. In the smooth case, Darboux theorem was a principal tool in Delzant's approach.
To analyze the local structures of the singular symplectic spaces in question, we find it necessary again
to make extensive use of the local normal form. The proof of Theorem ~TheoremIII~also requires the minimal coupling procedure of Sternberg \Sternberg~ and Weinstein \Weinstein.

In this paper, we only consider symplectic reductions for the zero level set. More tools will be needed to deal with other level sets. For a thorough review of more recent developments in reduction theory in the general case, see \OrtegaRatiu.

\newsec{Convexity of Image of Moment Map}

We begin with some notations, which will be used throughout the paper.
Let 
$$
1\ra A{\br i\over\ra} G{\br\pi\over\ra} T\ra1
$$ 
be an exact sequence of compact Lie groups where
$T=T^n$ is a torus. This induces the sequence of the dual of the Lie algebras
$0\la\ga^*{\br i^*\over\la}\gg^*{\br\pi^*\over\la}\gt^*\la 0$.
Let $M$ be a Hamiltonian $G$-space with a moment map $J$.
{\it The $A$-moment map $J_A:=i^*\circ J$ will be assumed proper throughout the paper.} Put
$$
X=J_A^{-1}(0)/A.
$$
This is a reduced $T$-space equipped with a $T$-moment map $J_T$ induced by $J$. 

Consider a fixed but arbitrary point $\tilde p\in J_A^{-1}(0)$ and let
$p=A\cdot\tilde p$. Let $V$ be the symplectic slice to the orbit $G\cdot\tilde p$, i.e.
$$
V=T_{\tilde p}(G\cdot\tilde p)^\omega/T_{\tilde p}(G\cdot\tilde p),
$$
the fiber at $\tilde p$ of the symplectic normal bundle of $G\cdot\tilde p$
in $M$. The symplectic form $\omega$ at the point $\tilde p$
induces a symplectic bilinear form on the vector space $V$
which we denote by $\omega_V$. Introduce the notations, all of which
depending on $\tilde p$,
$$\eqalign{
H&=Stab_G\tilde p\cr
K&=H\cap A\cr
B&=Image~of~H~under~(G\ra T)\cr
N&=T/B\cr
Q&=G/H\cr
L&=A/K.
}$$
Denote by $M_{(K)}$ the subset of points in $M$
whose stabilizers are conjugate to $K$ in $G$.
At the tangent space level,
we shall always denote by lower case gothic letters the corresponding
tangent spaces at the identity or identity cosets.
So we have the diagrams of exact sequences:
\eqn\MainDiagram{
\matrix{
 & & 0 & & 0 & & 0 & & &\cr
 & & \ua & & \ua & & \ua & & &\cr
0 & \ra &\gl&\ra&\gq&\ra&\gn&\ra&0\cr
 & & \ua & & \ua & & \ua & & &\cr
0 & \ra &\ga&\ra&\gg&\ra&\gt&\ra&0\cr
 & & \ua & & \ua & & \ua & & &\cr
0 & \ra &\gk&\ra&\gh&\ra&\gb&\ra&0\cr
 & & \ua & & \ua & & \ua & & &\cr
 & & 0 & & 0 & & 0 & & &\cr
}}
It is easy to see that $K$ is a normal subgroup of $H$, and that $1\ra B\ra T\ra N\ra 1$ is
an exact sequence of groups. Put
$$
D=Ker(G\ra N)
$$
where $G\ra N$ is the composition $G\ra T\ra N$. Then we have 
$$\eqalign{
D&=A\cdot H=H\cdot A\cr
D/H&\cong A/K\cr
H/K&\cong D/A.
}$$

\lemma{
$G$ contains a central torus $N'$ such that
$N'\ra N$ is a finite cover under $G\ra N$.}
\thmlab\FiniteCoverN
\proof
Since $N$ is connected, we may as well assume that $G$ is connected without loss of generality. Thus
$G$ has the shape $Z\times G_{ss}$ where $Z$ is the identity component of the center of $G$ and $G_{ss}$ the semi-simple part of $G$. It follows that $Z$ surjects onto $N$. Thus we have reduced our question to the case when $G$ is a torus. Thus we may as well write $Z=(S^1)^p=\R^p/\Z^p$ and 
$N=(S^1)^n=\R^n/\Z^n$. At the Lie algebra level we have a linear map $\pi:\R^p\thra\R^n$
where $\pi(\Z^p)=:\Pi\subset\Z^n$. Let $L=Ker(\Z^p\ra\Pi)$ and fix the standard inner product on $\R^p$. Then we have an exact sequence
$L_\R\hra \R^p=L_\R\oplus L^\perp_\R\thra\R^n$ where the last map is $\pi$ which
maps $L_\R^\perp$ isomorphically onto $\R^n$. Put $N'=L^\perp_\R/L^\perp$. Then
$N'\subset\R^p/\Z^p$ canonically. Since $L^\perp_\R\cong\R^n$, it follows that
$L^\perp\hra\Z^n$ is torsion. Hence $N'\ra N$ is a finite cover under $Z\ra N$. $\Box$


Obviously $H$ acts on $V$ linearly and symplectically.
It has a $H$-moment map $\Phi_V$ such that
$$
\xi\circ\Phi_V(v)=\half\omega_V(\xi_V\cdot v,v),~~~\xi\in\gh.
$$
Here $\xi_V$ is the operator on $V$ representing $\xi$.

\subsec{Local normal form for half-reduced space}

\lemma{(Local Normal Form) A neighborhood of the orbit $T\cdot p$
in $X$ is $T$-equivariantly symplectomorphic to a neighborhood
of the zero section of $Y_0=T\times_B(\gn^*\times\Psi^{-1}_V(0)/K)$
with the $T$-moment map $\tilde J_{T,v_p}:Y_0\ra\gt^*$ given by
$$
\tilde J_{T,v_p}([g,\eta,v])=\eta+\Phi_V(v)+v_p.
$$
Here $v_p=J_T(p)$, $\Psi_V:=j^*\circ\Phi_V$ is the induced moment map on $V$
for the subgroup $K{\br j\over\subset} H$.}
\thmlab\LocalNormalForm
\proof
The vector $v_p\in\gt^*\subset\gg^*$ is clearly $G$-invariant. By Proposition 2.5 \SL,
a neighborhood of the $G$-orbit $G\cdot\tilde p$ in $M$
is $G$-equivariantly symplectomorphic
to a neighborhood of the zero section of the vector bundle over $G/H$:
$$
Y=G\times_H(\gq^*\times V)
$$
with the moment map $\tilde J_{G,v_p}:Y\ra\gg^*$
$$
\tilde J_{G,v_p}([g,\eta,v])=Ad^*(g)(\eta+\Phi_V(v))+v_p.
$$
Here we have used a $G$-invariant inner product on $\gg^*$ to
make the identification $\gg^*=\gh^*\oplus\gq^*$. Since $i^*:\gg^*\ra\ga^*$
is a $G$-module homomorphism, we have
$i^*Ad^*(g)(\eta+\Phi_V(v))=Ad^*(g)(i^*\eta+i^*\Phi_V(v))$.
Since $Ad^*(g)$ is invertible and $i^*\eta\in\gl^*$, $i^*\Phi_V(v)\in\gk^*$,
it follows from \MainDiagram~ that $i^*\circ\tilde J_{G,v_p}([g,\eta,v])=0$
iff $\eta\in\gn^*, j^*\circ\Phi_V(v)=0$ where $j=i|H$. In other words,
$$
(i^*\circ\tilde J_{G,v_p})^{-1}(0)=G\times_H(\gn^*\times\Psi_V^{-1}(0)).
$$
Now taking the $A$-orbit space of this zero set,
we get
$$
(i^*\circ\tilde J_{G,v_p})^{-1}(0)/A=
T\times_B(\gn^*\times\Psi^{-1}_V(0)/K)=:Y_0.
$$
The $G$-equivariant symplectomorphism above restricts
and descends to a $T$-equivariant symplectomorphism from
a neighborhood of the $T$-orbit $T\cdot p=(G\cdot\tilde p)/A$
in $X=J_A^{-1}(0)/A$ to a neighborhood of
the zero section of $Y_0$.

Now $\tilde J_{T,v_p}$ is the map
induced by $\tilde J_{G,v_p}$ on $Y_0$.
For $[g,\eta,v]\in Y_0$, we have $Ad^*(g)=1$ and $\Psi_V(v)=0$,
hence $\Phi_V(v)\in\gb^*$. So we have
$$
\tilde J_{T,v_p}([g,\eta,v])=\tilde J_{G,v_p}([g,\eta,v])
=\eta+\Phi_V(v)+v_p\in\gt^*.~~~~ \Box
$$

\remark{As shown in \SL, the symplectic structure on the local normal form $G\times_H((\gg/\gh)^*\times V)$ is given by the standard form on $T^*G=G\times\gg^*$ plus the induced form $\omega_V$ on $V$.}
\thmlab\SSOnLocalNormalForm

Let $W$ be the symplectic orthogonal, with respect to $\omega_V$, of the
$K$-fixed subspace $V^K$ in $V$, so that 
$$
V=W\oplus V^K
$$
as $K$-modules. This is also an $H$-module decomposition, since $K$ is a normal subgroup of $H$.
Put
$$
\Psi_W:=\Psi_V|W.
$$

\lemma{We have $\Psi_V^{-1}(0)=\Psi_W^{-1}(0)\times V^K$. Hence 
$Y_0=T\times_B(\gn^*\times\Psi_W^{-1}(0)/K\times V^K)$.}
\thmlab\ProductForm
\proof
Let $v=(w,v_0)\in W\times V^K=V$.
Suppose $\Psi_V(v)=0$. Then $\omega_V(j(\xi)v,v)=0$ $\forall\xi\in\gk$
where $j:\gk\hra\gh$. Since $w,v_0$ are orthogonal and $v_0\in V^K$, it
follows that $\omega_V(j(\xi)w,w)=0$ implying that $w\in \Psi_W^{-1}(0)$.
The converse is similar. $\Box$

Note that each of the three factors in parentheses are $B=H/K$ invariant:
$B$ acts trivially on $\gn^*$ and $H$ leaves
each $\Psi_W^{-1}(0),V^K$ invariant, hence the $H$ action
descends to a $B$ action on the $K$ orbit spaces.

\lemma{In the local normal form of $\tilde p\in M$, $M_{(K)}\cap J_A^{-1}(0)$ corresponds to a 
relative neighborhood of the zero section in $G\times_H(\gn^*\times0\times V^K)$.
$T\times_B(\gn^*\times 0\times V^K)$ is the unique
stratum containing $[e,0,0,0]$ in the stratified space $Y_0$.}
\thmlab\GTStratum
\proof
Recall that the local normal form at $\tilde p\in J^{-1}(0)$ is
$Y=G\times_H(\gq^*\times V)$. The same is true if $\tilde p\in J_A^{-1}(0)$, since $G$ acts trivially on $\gt^*$. The stratum $M_{(K)}$ of orbit type $K$
is locally $Y_{(K)}=G\times_H(\gq^*_{(K)}\times V_{(K)})$,
Recall that $J_A^{-1}(0)$ is locally $G\times_H(\gn^*\times\Psi_V^{-1}(0))$.
We claim that
$$
\gn^*\subset\gq^*_{(K)},~~~~V_{(K)}\cap\Psi_V^{-1}(0)=V^K.
$$
Since $T$ acts trivially on $\gt^*$,
so does $G$ via the homomorphism $G\ra T$.
In particular $K=H\cap A$ also
acts trivially on $\gn^*\subset\gt^*$.
This shows that $\gn^*\subset(\gq^*)^K\subset
\gq^*_{(K)}:=\cup_{g\in H} (\gq^*)^{g^{-1}Kg}$.
Now $V_{(K)}=\cup_{g\in H}V^{g^{-1}Kg}$. But $V^{g^{-1}Kg}=V^K$
because $H\cdot V^K=V^K$. This shows that $V_{(K)}=V^K$.
On the other hand, we have
$\Psi_V^{-1}(0)=\Psi_W^{-1}(0)\times V^K\supset V^K$. This proves the
equality above. It follows that
$M_{(K)}\cap J_A^{-1}(0)$ is locally
$$
G\times_H(\gn^*\times0\times V^K).
$$
Taking quotient of this by $A$, we see that the
local normal form of $(M_{(K)}\cap J_A^{-1}(0))/A$ is
$T\times_B(\gn^*\times 0\times V^K)$.
The former space is the unique stratum
in $X=J_A^{-1}(0)/A$ containing $p$ (Theorem 2.1 \SL). Our assertion
is nothing but the local version of this. $\Box$

\subsec{Proof of the convexity property}

For any closed subgroup $C\subset T$, the moment map $J_T:X\ra\gt^*$
for the $T$ action on $X=J_A^{-1}(0)/A$
induces a moment map for the $C$ action which we denote by $J_C$.

\lemma{For any $\eta\in\gc^*$, the level set $J_C^{-1}(\eta)$ is connected.}
\proof
Let $G'\subset G$ be the preimage of $C\subset T$ under the map $G\ra T$,
and let $J_{G'}$ be the induced $G'$-moment map. Since $A\subset G'$ induces a projection $r^*:{\gg'}^*\ra\ga^*$ and $J_A=r^*\circ J_{G'}$, it follows that  $J_{G'}$ is proper, because $J_A$ is.
Since $G$ acts trivially on $\gt^*$, we can view $\eta$ as a $G'$-invariant element of ${\gg'}^*$.
By a result of Kirwan \Kirwan, $J_{G'}^{-1}(\eta)$ is connected. It follows that 
$J_C^{-1}(\eta)=J_{G'}^{-1}(\eta)/A$ is connected. $\Box$

\theorem{The moment map $J_T:X\ra\gt^*$ has the following properties:
\item{i.} the level sets of $J_T$ are connected;
\item{ii.} the image of $J_T$ is convex;
\item{iii.} the image of $J_T$ is the convex hull of the image of the
fixed points of the $T$ action.}
\proof
Part i. follows from the preceding lemma.
The main point of the rest of the proof
is that the strategy of the proofs of Atiyah and Guillemin-Sternberg in the smooth case
carry over, but with two changes.
The fixed point set is now a stratified space
(rather than a manifold), and the local structure of
a fixed point is replaced by the local normal form in a stratified space.
We will use the identifications $T=\R^n/\Z^n$
and $\gt=\R^n\equiv\gt^*$.

Part ii. Let $\Delta=J_T(X)$. Since $X$ is compact, so is $\Delta$. For every line
$L=\{v_0+tv_1|t\in\R\}$ with rational direction i.e. $v_1\in\Q^n$,
we will show that $L\cap\Delta$ is connected. Consider the $(n-1)$
dimensional Lie subgroup $C$ with Lie algebra
$\{a\in\R^n|\bra v_1,a\ket=0\}$, and let $P:\R^n\ra\R^n/\R v_1$
be the projection. Then $J_C=P\circ J_T$ is a moment map for
the $C$ action on $X$. By the preceding lemma, $J_C^{-1}(Pv_0)$
is connected. So $J_T(J_C^{-1}(Pv_0))=P^{-1}(Pv_0)\cap J_T(X)=L\cap\Delta$
is connected. This shows that $\Delta$ is convex.

Part iii. By \SL, the fixed point set $X^T$ is closed subset
which is a disjoint union of closed
connected stratified symplectic subspaces $C_i$ of $X$.
We claim that there are only finitely many connected
components $C_i$. Assume the contrary.
Then there is an infinite set of points $p_i\in C_i$.
Let $p$ be a limit point of this set. By continuity of the $T$ action,
we have $p\in X^T$, hence $p\in C_i$ for some $i$.
By the local normal form, there is a neighborhood
of $T\cdot p$ in $X$ which is equivariantly symplectomorphic to
a neighborhood of the zero section of
$Y_0=T\times_B(\gn^*\times\Psi_V^{-1}(0)/K)$.
But since $T\cdot p=p$, it follows that $B=T$, $N=1$,
and $p$ corresponds to $0\in\Psi_V^{-1}(0)/K$.
Note that if $x\in\Psi_V^{-1}(0)/K$ is any $T$ fixed point
then it is connected to $0$ because $T$ acts linearly on $V$
and the points $\{tx|t\in\R\}$ are $T$-fixed.
This shows that the points $p_i$ sufficiently close to
$p$ must all be in the same connected component of fixed points,
contradicting that the $C_i$ are distinct connected components of
$X^T$. Hence $X^T$ is a union of finitely many components $C_1,..,C_N$.

Since the $T$-action on each stratified space $C_i$
is trivial, $J_T|C_i$ must be constant. Put $\eta_i=J_T(C_i)$.
By convexity of $\Delta$, the convex hull $\Delta'=conv(\eta_1,..,\eta_N)$
is a subset of $\Delta$. Suppose that $\xi\in\Delta\backslash\Delta'$.
We can choose $\xi$ so that the $\xi-\eta_i$ all lie in
the same half space bounded by a hyperplane in $\R^n$.
Let $\chi$ be a normal vector to the hyperplane so that
$$
\bra\xi,\chi\ket>\bra\eta_i,\chi\ket~~\forall i.
$$
Choose the hyperplane so that $\chi$ is generic i.e. the components of
the vector $\chi$ are independent over $\Q$. Then
the one-parameter subgroup $T'=\{exp~t\chi|t\in\R\}$
is a dense subgroup in $T^n$. Note that the zeros
of the vector field $V_\chi$ on $X$ are $T'$ fixed points,
hence $T$ fixed points (by density).

Let $p$ be a point where the function
$\bra J_T,\chi\ket$ on $X$ attains a maximum, say at $p$. Then in the stratum containing $p$, we have
$\bra dJ_T,\chi\ket=0$ at $p$. By the moment map condition,
it follows that the vector field $V_\chi$ vanishes at $p$,
hence $p$ is a fixed point, which means that $J_T(p)=\eta_i$
for some $i$. This implies that $\bra \eta_i,\chi\ket
\geq\langle J_T(x),\chi\rangle$ for all $x\in X$,
contradicting the inequalities above.
This shows that $\Delta=\Delta'$. $\Box$

\newsec{Rationality of Moment Polytope}

In this section, we prove the first assertion in Theorem \TheoremII.
The notations introduced in the last section, such as \MainDiagram, will remain in force here.
Thus $X$ is the reduced Delzant $T$-space obtained by reducing
a connected Hamiltonian $G$-manifold $(M,\omega)$ with respect
to a normal subgroup $A\subset G$. First we prove that the moment polytope classifies the
$T$-orbits in $X$.

\subsec{Orbit theorem}

\lemma{$B$ acts effectively on $E:=\Psi^{-1}_W(0)/K\times V^K$
and $2dim~B=dim~E$.}
\proof
Note that $B$ acts trivially on $\gn^*$.
Suppose $C\subset B$ is a subgroup that acts trivially on the open dense stratum of $E$.
Then we have $Y_0=T\times_B (\gn^*\times E)\cong T/C\times_{B/C}(\gn^*\times E)$.
Since $T$ is abelian, this means that $C$ acts trivially on this
fiber product. We claim that $C=1$, which shows that $B$ acts effectively on $E$.

When a Lie group acts effectively on a space in a Hamiltonian fashion, then any subgroup of $T$ that
fixes a nonempty open subset must be the trivial group.
By assumption $T$ acts effectively on the open dense stratum of $X$.
It follows that $T$ acts effectively on a neighborhood of the zero section of $Y_0$.
Since $C$ acts trivially on $Y_0$, it follows that $C=1$.

Finally the dimension assertion follows from that $2dim~T=dim~X=dim~Y_0=dim~T+dim~\gn+dim~E-dim~B$,
and that $dim~T=dim~\gn+dim~B$. $\Box$

\corollary{The generic $B$ orbits in $\Psi_W^{-1}(0)/K$
and in $V^K$ have dimensions exactly half the respective dimensions of
those spaces.}
\proof
Call the respective dimensions of those symplectic spaces $2a,2b$,
and consider a generic orbit $B\cdot(p,q)$ in $E$. By the preceding lemma,
this orbit has dimension $dim~B=\half dim~E=a+b$. Since
$B\cdot(p,q)\subset B\cdot p\times B\cdot q$. Since the dimension of
an orbit of a symplectic action on a stratified space cannot exceed half the dimension of the space,
it follows that $dim~B\cdot p\leq a$ and $dim~B\cdot q\leq b$.
But these must be equalities in order that $dim~B\cdot(p,q)=a+b$. $\Box$


\lemma{Let $C=Ker(B\ra Sp(V^K))$. Then $2dim~B/C=dim~V^K$. 
Moreover, $V^K$ has no nonzero $B$-fixed point.}
\thmlab\VKFixedPoint
\proof
%
Since $B/C$ acts effectively and symplectically on $V^K$,
the generic orbit has dimension $dim~B/C\leq\half dim~V^K$.
But this must be an equality since $2dim~B=dim~E$.

Now since $B/C$ is a torus acting linearly and effectively on a linear space $V^K$ with $2dim~B/C=dim~V^K$, the only fixed point is $0$. $\Box$

\corollary{Let $S$ be the stratum containing $p$ in $X$.
Then $dim~S^B=2dim~T/B$.}
\proof
By Lemma \LocalNormalForm, a neighborhood of $p$ in $X$ is $T$-equivariantly
symplectomorphic to a neighborhood of $[e,0,0,0]\in T\times_B(\gn^*\times
\Psi_W^{-1}(0)/K\times V^K)$. By Lemma \GTStratum, a neighborhood
of $p$ in $S$ is mapped into the stratum $E:=T\times_B(\gn^*\times 0\times V^K)$ in $Y_0$. Thus
we have $dim~S^B=dim~E^B$.
Since $T$ is abelian, we have
$E^B=T\times_B(\gn^*\times 0\times(V^K)^B)$.
By the preceding corollary, $(V^K)^B=0$.
This shows that $E^B\cong T/B\times\gn^*$,
which has dimension $2dim~T/B$ because $\gn=\gt/\gb$. $\Box$

\corollary{$\Psi_W^{-1}(0)/K$ has no nonzero $B$-fixed point.}
\proof
If $x\in\Psi_W^{-1}(0)/K$ is a fixed point, then so is the set $\R\cdot x$.
Suppose $x\neq 0$, which we may assume to be close to $0$.
The point $[e,0,0,0]\in T\times_B(\gn^*\times x\times 0)$,
correspond to some point $q\in X$.
Let $S$ be the stratum containing $q$ in $X$. Then
$dim~S^{B}=2dim~T/B$ by the preceding corollary.
But as we move along $\R\cdot x$ (in the stratum containing $[e,0,x,0]$ in the local model),
the $B$-fixed point set around this point
in $T\times_B(\gn^*\times\R\cdot x\times 0)$
would have dimension at least $2dim~T/B+1$, a contradiction. $\Box$

\corollary{The $T$-fixed points of $X$ are isolated.}
\proof
Note that at a $T$-fixed point $p$, we have $B=T, H=G, K=A$, 
so that the local normal form is a neighborhood of
$(0,0)\in\Psi_W^{-1}(0)/K\times V^K$.
By the preceding lemma, $V^K$ has no nonzero $B$-fixed point.
By the preceding corollary, the same holds for $\Psi_W^{-1}(0)/K$.
So $(0,0)$ is the only $B$-fixed point. $\Box$

\lemma{For $v\in V^K$, if $\Phi_V(v)=0$ then  $v=0$.}
\thmlab\PhiVanishLemma
\proof
Let $C=Ker(B\ra Sp(V^K))$ as before. We claim that $V^K\ra(\gb/\gc)^*\subset\gb^*$,
$v\mapsto\Phi_V(v)$ is the moment map for the linear $B$-action on $V^K$ sending $0\mapsto 0$.
This suffices, for then the moment map must be of the form
$v\mapsto \Phi_V(v)=(\sum_{j=1}^{dim~B/C}w_{ij}|z_j(v)|^2)$
where the $z_j$ are a choice of the linear complex coordinates of $V^K$ and $(w_{ij})$ is some matrix.
Since $B/C$ acts effectively, this matrix has full rank and we know that $2dim~B/C=dim_\R~V^K$,
it follows that $\Phi_V(v)=0\Lra |z_j(v)|^2=0$ for all $j$ hence $v=0$.

We now prove the claim. By definition $\Phi_V:V\ra\gh^*$
is the unique $H$-moment map for $V$ with $\Phi_V(0)=0$.
Since $\Psi_V^{-1}(0)/K$ is a connected stratified Hamiltonian $B$-space,
$\Phi_V$ induces a moment map by restricting to the zero level set and taking $K$-orbits.
By Lemma \ProductForm, $\Psi_V^{-1}(0)/K=\Psi_W^{-1}(0)/K\times V^K$,
we can further restrict the moment map to the stratum $0\times V^K$, and $\Phi_V$
induces a $B=H/K$-moment map.
Since $C\subset B$ acts trivially, $\Phi_V$ further induces
a $B/C$-moment map on $V^K$. This completes the proof. $\Box$

\theorem{(Orbit Theorem, cf. p21 \Guillemin)
The map $J_T:X\ra\Delta:=J_T(X)$ descends to
a homeomorphism $X/T\ra\Delta$.}
\proof
Since $T$ is compact and $J_T$ is a $T$-equivariant continuous map of
compact stratified spaces, if the quotient map $X/T\ra\Delta$
is a bijection it is automatically a homemorphism.
Given $p\in X$, put $v_p=J_T(p)$. By $T$-equivariance of $J_T$,
$J_T^{-1}(v_p)$ contains at least one $T$-orbit, and we want to show that $J_T^{-1}(v_p)$ contains no more than one orbit. 
Consider the local normal form of $p$, Lemma \LocalNormalForm, given by
$Y_0=T\times_B(\gn^*\times\Psi_V^{-1}(0)/K)$ with moment map
$\tilde J_{T,v_p}([g,\eta,v])=\eta+\Phi_V(v)+v_p$.
Now $\tilde J_{T,v_p}([g,\eta,v])=v_p$ iff $\eta+\Phi_V(v)=0$
iff $\eta=0=\Phi_V(v)$ because $\eta\in\gn^*$, $\Phi_V(v)\in\gb^*$. We know that $T\cdot p$ must be contained in the stratum of $p$.
The stratum of $p$ in $Y_0$ is $T\times_B(\gn^*\times 0\times V^K)$ by Lemma \GTStratum. So for a point $[g,\eta,v]$ to be in the orbit of $p$, we may assume that $v$ lies in $V^K$. By the preceding lemma, we have $v=0$. This shows that the $T$-orbit of $p$ contained in that stratum and in $J_T^{-1}(v_p)$ must be $T\times_B(0\times 0\times 0)$. This shows that each stratum of $X$ contains at most one $T$-orbit in $J_T^{-1}(v_p)$. But since $J_T^{-1}(v_p)$ is connected by Theorem \TheoremI i, there can't be more than one $T$-orbit altogether. $\Box$

\corollary{If $F\subset\Delta$ is any connected subset,
then $J^{-1}(F)$ is connected.}
\proof
By the preceding lemma $J^{-1}(F)/T$ is connected.
Since $T$ is connected, it follows that $J^{-1}(F)$ is connected.
$\Box$

\subsec{Rationality}

{\it Notation.} For $p\in X$, $v_p=J_T(p)$ lies in the relative interior of a
unique face of the moment polytope $\Delta=J_T(X)$. We denote that interior
of that face by $F$. We shall refer to $F$ as the face containing $v_p$. 
Note that the interior $\Delta^\circ$ of $\Delta$ is the largest face.
Again, the reader is reminded that the notations $B,D,H,K,L,Q,N,V,W,F$
are all associated with the given point $p\in X$. If $p'\in X$
is a second point, we denote those associated objects by $B',D',...$.

\lemma{Let $S$ be the stratum containing $p$ in $X$.
If $S_0$ is the connected component containing $p$ in $S^B$
then $J_T(S_0)$ is a connected open subset of the affine subspace $J_T(p)+\gn^*$ in $\gt^*$.}
\thmlab\DimensionStratum
\proof
By Lemma \LocalNormalForm,
in some neighborhood $U_0\ni p$ in $X$, we know that
$J_T:U_0\cap S_0\ra\gt^*$
is represented by the projection $T\times_B(\gn^*\times 0\times 0)\ra\gn^*
\subset\gt^*$ plus $v_p=J_T(p)$.
Thus $J_T(U_0\cap S_0)=J_T(p)+O_0$ for some open
neighborhood $O_0\ni 0$ in $\gn^*$.  By Lemma \GTStratum, every point
$q$ in a relative neighborhood of $p$ in $S_0$ has the same $B,\gn$, etc.
Since $S_0$ is connected, we can cover it with small open sets
and repeat the argument above in each open set. In the end, we see that $J_T(S_0)$
is a union of sets of the form $J_T(q)+O$, where $q\in S_0$ with $J_T(q)-J_T(p)\in\gn^*$ and $O$ is a relative neighborhood of $0$ in $\gn^*$. This shows that $J_T(S_0)$ is an open subset of $J_T(p)+\gn^*$. Since $S_0$ is connected, so is $J_T(S_0)$. $\Box$.

\corollary{For any $p\in X$, the face $F$ containing $v_p$ also contains an open
neighborhood of $v_p$ in $v_p+\gn^*$.}
\proof
Since $J_T(S_0)\subset\Delta$, there is a neighborhood
of $v_p$ in $v_p+\gn^*$ which lies in $\Delta$. Since $\Delta$
is convex by Theorem \TheoremI ii, this neighborhood must lie in $F$. $\Box$

\lemma{
Let $F$ be the interior of a given face of $\Delta$.
Let $p\in J_T^{-1}(F)$
be a point such that $dim~\gn$ is the largest possible
as $p$ varies over $J_T^{-1}(F)$
(remember that $\gn$ is associated to $p$). Then
in the local normal form, a neighborhood of $p$ in $J_T^{-1}(F)$ corresponds to the stratum $T\times_B(\gn^*\times 0)$. Moreover $F\subset v_p+\gn^*$.}
\proof
In a small neighborhood $U_0$ of $p$ in $X$,
we may regard $J_T^{-1}(F)\cap U_0$ as a subset of
$Y_0=T\times_B(\gn^*\times\Psi_V^{-1}(0)/K)$ by Lemma \LocalNormalForm. We claim that
$p'=[t,\eta,v]\in J_T^{-1}(F)\cap U_0\Lra v=0$.
Otherwise the stabilizer of $p'$ in $T$ would be a proper
subgroup $B'\varsubsetneq B$ by Lemma \VKFixedPoint~ and one
of its corollaries. The local normal form of $p'$ would
then take the shape $T\times_{B'}({\gn'}^*\times\Psi_{V'}^{-1}(0)/K')$. Since $B'\varsubsetneq B$,
this would mean that $dim~\gn'>dim~\gn$ and that $F$ contains
a neighborhood of $v_{p'}$ in $v_{p'}+\gn'$,
contradicting the maximality of $dim~\gn$. It follows that
$p'=[t,\eta,0]$. So we have shown that
$J_T^{-1}(F)\cap U_0\subset T\times_B(\gn^*\times 0)$.
Applying $J_T$ to both sides we get that
$F\subset v_p+\gn^*$. $\Box$

{\it Warning.} The argument above only establishes that
$F\subset v_p+\gn^*$ for some
$p\in J_T^{-1}(F)$. It does not show that $F\subset v_p+\gn^*$
for every $p\in J_T^{-1}(F)$. For this, we need more work below. 

\theorem{$J_T(X)=\Delta$ is a rational polytope in $\gt^*$.}
\proof
Each facet of $\Delta$ lies in an affine hyperplane of the shape $v_p+\gn^*$, which has a normal vector in the line $\gb=Lie(B)$. It follows that $\Delta$ is rational. $\Box$

\subsec{Face theorem}

\theorem{Let $F$ be the interior of a given face in $\Delta$.
For any $p\in J_T^{-1}(F)$ we have $F\subset v_p+\gn^*$, hence $dim~F=dim~\gn^*$.}
\proof
The dimension assertion follows from the fact that
$F$ contains a neighborhood of $v_p$ in $v_p+\gn^*$
and the inclusion assertion, which we now show.
Suppose the contrary, which means that $dim~F>dim~\gn^*$.
By the preceding theorem $F$ is rational, hence there is a
{\it rational} affine line $\ell$ transverse to $\gn^*$ in $F$.
Note that $\ell\cap(v_p+\gn^*)=v_p$.
Since $\gt':=(\ell-v_p)^\perp\subset\gt$ is a rational subspace,
it is the Lie algebra of a codimension one closed subtorus $T'\subset T$.

Let $G'$ be the preimage of $T'$ in $G$ and $J_{G'}$ the $G'$-moment map
induced by $J$. Then we have $J_A=r^*\circ J_{G'}$, $r^*:{\gg'}^*\ra\ga^*$ is the projection
induced by $A\subset G'$. Since $J_A$ is proper, so is $J_{G'}$.
The reduced space $X':=J_{G'}^{-1}(w)/G'$
is a stratified symplectic space with
a unique open connected dense stratum (Theorem 5.9 \SL). Here 
$w$ is the image under ${\gt'}^*\hra{\gg'}^*$ of $j^*v_p\in{\gt'}^*$,  where $j:\gt'\hra\gt$.
Note that $w$ is a $G'$-invariant vector because $G'$ acts trivially on $\gt'$ via $G'\ra T'$.
By section 4 \SL, we have a moment map $J_{T'}:X\ra{\gt'}^*$, and
the space $X'=J_{T'}^{-1}(j^*v_p)/T'$ has an induced $T/T'$ action with a moment map
$J_{T/T'}:X'\ra{\gt'}^\perp=\ell-v_p$. Since $dim~\gt'=dim~\gt-1$ and $dim~X=2dim~\gt$
and $T$ acts effectively, it follows that $dim~X'=2$. The image $J_{T/T'}(X')$ is $\ell\cap\Delta=\ell\cap F_p$, a line segment. By the Orbit Theorem, $X'/(T/T')$ is
homeomorphic to this line segment.

Now the complement of the point $v_p$ in $\ell$
is a disjoint union of open half lines
$\ell_1,\ell_2$. In particular, the $J^{-1}_{T/T'}(\ell_i)$
are two non-empty open subsets of $X'$. We also have
$$
J_{T/T'}^{-1}(v_p)=J_T^{-1}(v_p)/T'=T\cdot p/T'\cong T/BT'.
$$
Since $T$ is connected, so is $T/BT'$. It is a point
because the tangent space at $e$ is $\gt/(\gb+\gt')$, which is zero.

So $X'$ is a two dimensional space
which is a disjoint union of two nonempty open sets
$J_{T/T'}^{-1}(\ell_i)$ and a single
point $J_{T/T'}^{-1}(v_p)$. But $X'$ has a connected open
dense stratum, which is a contradiction. $\Box$

\corollary{$dim~J_T^{-1}(v_p)=dim~F$.}
\proof
By the Orbit Theorem, we have
$J_T^{-1}(v_p)=T\cdot p\cong T/B=N$, which has dimension $dim~F$
by the preceding theorem. $\Box$

\corollary{$p\in X^T$ iff $F$ is a vertex of $\Delta$.}

\corollary{$dim~J_T^{-1}(F)=2dim~F=2dim~T/B$.}

\corollary{Let $p\in J_T^{-1}(F)$, $S$ be the stratum containing $p$ in $X$,
and $S_0$ the connected component containing $p$ in $S^{B}$. Then $J_T^{-1}(F)=S_0$.}
\proof
By Lemma \DimensionStratum, $J_T(S_0)$ is a connected open set
in the (relative) interior of $\Delta\cap(v_p+\gn^*)$. By the Face Theorem,
$F\subset v_p+\gn^*$, it follows that
$F$ is the interior of $\Delta\cap(v_p+\gn^*)$.
Hence $J_T(S_0)\subset F$, i.e. $S_0\subset J_T^{-1}(F)$.
Note that $S_0$, being a connected component of a fixed point set $S^B$ in
a stratum, is a manifold of dimension $2dim~T/B$ (a corollary to Lemma \VKFixedPoint). Therefore $J_T^{-1}(F)$ is a union of such manifolds $S_0'$,
one for each $p'\in J_T^{-1}(F)$, and any two are either disjoint or equal. 
But since $J_F^{-1}(F)$ is connected and has the same dimension as each $S_0'$,
these manifolds $S_0'$ must all coincide. Thus $J_T^{-1}(F)=S_0$. $\Box$

\corollary{The points in $J_T^{-1}(F)$ have the same stabilizer
group in $T$.}
\thmlab\SameStabilizerLemma
\proof
Near each point $p\in J_T^{-1}(F)$, $S_0$ looks locally like $T\times_B(\gn^*\times 0)$
where every point has stabilizer $B$. So there is a neighborhood of $p$ in $S_0$
having the same stabilizer. Thus the set of
points in $J_T^{-1}(F)$ having a given stabilizer group
is open in $S_0$. Since $S_0$ is connected, there cannot
be two different stabilizer groups. $\Box$

\corollary{$J_T^{-1}(\Delta^\circ)$ is the open dense stratum in $X$ and $T$ acts freely there.}
\proof
Since $X$ is a reduced Delzant $T$-space, $T$ acts effectively on the
unique open dense stratum in $X$. Since the strata in $X$ are labeled by stabilizer
subgroups $B$ in $T$, and each $J_T^{-1}(F)$ has dimension $2dim~T/B$,
it follows that the open stratum must correspond to a $B$ which is finite
in order that the stratum has dimension $2dim~T=dim~X$. But by uniqueness of the open stratum,
and the effectiveness of the $T$-action, $B=1$, and $J_T^{-1}(\Delta^\circ)$ is the
only stratum with this property. $\Box$

\corollary{Let $\pi:J_A^{-1}(0)\ra X$ be the projection.
The points in $\pi^{-1}J_T^{-1}(F)=J_G^{-1}(0\times F)$
have the same stabilizer type in $G$.}
\proof
This follows from the next lemma and the fact that $F$ is
connected. $\Box$

\lemma{In the local normal form of a point $\tilde p\in J^{-1}(0\times F)$,
we have $J^{-1}(0\times F)\cap U_0\subset G\times_H(\gn^*\times 0)$.
In particular if $O$ is a small open set in $F$,
then the points in $J^{-1}(0\times O)$ have the same stabilizer type in $G$.}
\thmlab\StabilizerType
\proof
Again we represent a neighborhood of $\tilde p$ in $M$ by
its local normal form. Then the composition map
$J_A^{-1}(0){\br\pi\over\ra} X{\br J_T\over\ra}\gt^*$, which is a restriction of $J:M\ra\gg^*\equiv\ga^*\oplus\gt^*$, is represented by 
$$
G\times_H(\gn^*\times\Psi_V^{-1}(0))\ra\gt^*,~~~
[g,\eta,v]\mapsto Ad^*(g)(\eta+\Phi_V(v))+v_p=\eta+\Phi_V(v)+v_p
$$
by Lemma \LocalNormalForm. Since $F\subset\gn^*+v_p$,  by the Face Theorem, it follows that
$[g,\eta,v]\in J^{-1}(0\times F)\Lra \Phi_V(v)=0$. This shows that for a small open set $O\subset F$,
$J^{-1}(0\times O)$ is (corresponding to) a subset of
$G\times_H(\gn^*\times\Phi_V^{-1}(0))$. But the image of $J^{-1}(0\times O)$ under
$\pi$ is $J_T^{-1}(O)$, where every point has the same stabilizer $B$ in $T$, by Corollary
\SameStabilizerLemma. This means $J_T^{-1}(O)$ must be (corresponding to)
a subset of $T\times_B(\gn^*\times 0)$. This implies that for $[g,\eta,v]\in J^{-1}(0\times F)$
we must have $v=0$. In other words, $J^{-1}(0\times O)$ is a subset of $G\times_H(\gn^*\times 0)$.

Now since $H$ acts trivially on $(\gg/\gh)^*\hla\gn^*$,
it is easy to check that each point in $G\times_H(\gn^*\times 0)$
has stabilizer which is conjugate to $H$ in $G$. $\Box$

\corollary{The restriction of $\omega$ to $J^{-1}(0\times F)\subset M$
has constant rank.}
\proof
By the preceding lemma, $\omega$ restricted to $J^{-1}(0\times F)$
is locally equivalent to the symplectic form on the local
normal form restricted to $G\times_H(\gn^*\times 0)$.
Since $\gn$ remains the same for all $p\in J_T^{-1}(F)$, by Corollary \SameStabilizerLemma, 
it follows that the symplectic form has constant rank. $\Box$

\newsec{Realizing Reduced Delzant Spaces}

In this section, we shall prove the second part of Theorem \TheoremII.
Thus given an $n$ dimensional
rational polytope $\Delta$ in ${\R^n}^*$,
we will construct a reduced Delzant $T$-space $X_\Delta$,
whose moment polytope coincides with $\Delta$. Here $T=\R^n/\Z^n$.
In the next subsection, we will show that $X_\Delta$ also has the
structure of a complete toric variety. Both constructions are generalizations of the
construction for the case when $\Delta$ is regular or simplicial (see Guillemin's book \Guillemin).

There is one subtle point that requires clarification. The notion of a rational polytope, Definition \DefRationality, is in terms of a choice of a torus $T$. If $T=\R^n/\Pi$, then $\gt^*=Lie(T)^*={\R^n}^*$.
A polytope $\Delta$ in $\gt^*$ is rational iff each facet has a normal vector $u$ which lies in $\Pi$.
This is equivalent to requiring that a normal vector $u$ lies in $\Z^n$.
In particular the same set $\Delta\subset {\R^n}^*$ can be rational with respect to many different tori $T$ with the same Lie algebra. Moreover, the same set $\Delta$ can be the moment polytope of many $T$-spaces for different $T$ with the same Lie algebra. These different $T$-spaces need not even be homeomorphic (see below). This does not contradict Theorem \TheoremIII, which is a statement about $T$-spaces for a single chosen $T$.

\subsec{From rational polytopes to reduced spaces}

First let's label the codimension one faces of $\Delta$:
$1,..,p$. The $i$th one
lies in a unique hyperplane $\bra x,u_i\ket=\lambda_i$
where $u_i$ is the unique primitive vector in $\Z^n$
which is an inward pointing normal of the hyperplane. In other words,
$x\in\Delta$ iff
$$
\bra x,u_i\ket\geq\lambda_i,~~~\forall i.
$$
Define $\pi:\Z^p\ra\Z^n$, $e_i\mapsto u_i$, so that there is an exact
sequence
$$
0\ra L{\br\iota\over\hra}\Z^p{\br\pi\over\ra}\pi(\Z^p)=:\Pi\ra0
$$
where $L\subset\Z^p$ is the set of vectors
$l$ such that $\sum_il_iu_i=0$. Note that since $\Delta$ is
assumed $n$ dimensional $\pi:\R^p\ra\R^n$ is surjective.
But $\pi:\Z^p\ra\Z^n$ is surjective iff the $u_i$ generates $\Z^n$. Put
$$
A:=L_\R/L,~~~ T:=\R^n/\Pi,~~~G:=(S^1)^p=\R^p/\Z^p.
$$
Then we have an exact sequence of Lie groups $1\ra A\ra G\ra T\ra1$, and
canonical identifications $\ga=L_\R$,
$\gt=\R^n$, $\gg=\R^p$.
Let $\C^p$ be given the standard symplectic form,
and let $G$ act on $\C^p$ by the usual coordinate-wise scaling.
A moment map $J_G:\C^p\ra\gg^*={\R^p}^*$ is given by
$$
J_G(z)=\half(|z_1|^2,...,|z_p|^2)+\lambda
=\half\sum_i(|z_i|^2+2\lambda_i)e_i^*\in{\R^p}^*_\geq+\lambda.
$$

\lemma{$J_G(\C^p)={\R^p}^*_\geq+\lambda$
and ${\iota^*}^{-1}(0)=L^\perp$ in ${\R^p}^*$.}
\proof
Straightforward. $\Box$

\lemma{$\pi^*\Delta=({\R^p}^*_\geq+\lambda)\cap L^\perp
=({\R^p}^*_\geq+\lambda)\cap{\iota^*}^{-1}(0)$ in ${\R^p}^*$.}
\thmlab\MomentImage
\proof
Note that $L^\perp=Ker(\iota^*:{\R^p}^*\ra L^*_\R)$.
First $\pi^*\Delta\subset Im~\pi^*=Ker~\iota^*$ implies
$\iota^*$ kills $\pi^*\Delta$, which means that $l\cdot\pi^*\Delta=0$,
$\forall l\in L$. Thus $\pi^*\Delta\subset L^\perp$. We have
$$
x\in\Delta\LRa \bra x,\pi(e_i)\ket\geq\lambda_i,~\forall i
\LRa\bra\pi^*(x)-\lambda,e_i\ket\geq0,~\forall i
\LRa \pi^*(x)-\lambda\in{\R^p}^*_\geq.~~~~\Box
$$

\lemma{$J_A^{-1}(0)=(\iota^*\circ J_G)^{-1}(0)=J_G^{-1}(\pi^*\Delta)$.}
\thmlab\ZeroSet
\proof
We have
$$
J_G^{-1}{\iota^*}^{-1}(0)=J_G^{-1}(L^\perp)
=J_G^{-1}(L^\perp\cap({\R^p}^*_\geq+\lambda))
=J_G^{-1}(\pi^*\Delta).~~~~\Box
$$

\lemma{$J_A:=\iota^*\circ J_G:\C^p\ra\ga^*=L_\R^*$ is proper.}
\proof
Let $C\subset\ga^*$ be a compact set. Since $C$ is closed
so is $J_A^{-1}(C)$. So it suffices to show that this set is bounded.
Since any linear
projection can be represented topologically by orthogonal
projection, we can always find a closed ball $B\subset{\R^p}^*$
such that $C\subset\iota^*(B)$. We will show that $J_A^{-1}(\iota^*(B))$
is bounded. We have
$$
J_A^{-1}(\iota^*(B))=\cup_{b\in B}J_A^{-1}(\iota^*(b))
=\cup_{b\in B}J_G^{-1}(\pi^*\Delta_b)
=J_G^{-1}\pi^*(\cup_{b\in B}\Delta_b)
$$
where $\Delta_b$ is the set defined by
$\bra x,u_i\ket\geq\lambda_i-b_i$. Here we have used that
$J_A^{-1}(\iota^*(b))=J_G^{-1}(\pi^*\Delta_b)$ by a computation
similar to the preceding lemma.

Note that since $\Delta_0=\Delta$ is assume a convex polytope,
each $\Delta_b$ (possibly empty)
remains bounded because it can be obtained from
$\Delta_0$ simply by parallel translating bounding planes.
Since $b$ is varying over a bounded set, the union of polytopes
$\Delta_b$ remains bounded.
Since $\pi^*$ is linear, it follows that the
$\pi^*(\cup_{b\in B}\Delta_b)$ is also bounded.
Finally since $J_G$ is proper, the inverse image of
a bounded set is bounded. This completes the proof. $\Box$

\lemma{The zero set $J_A^{-1}(0)$ has dimension $p+n$
and $A$ acts effectively on the set. Hence $X_\Delta:=J_A^{-1}(0)/A$
has dimension $2n$. Moreover the induced $T$ action on $X_\Delta$
is effective.}
\proof
Note that $dim~A=p-n$, so that $dim~X_\Delta=2n$ follows from
our first assertion, which we now prove. 

Since $dim~{\iota^*}^{-1}(0)=n$ and each fiber $J_G^{-1}(y)$
has dimension at most $p$, it follows that $J_A^{-1}(0)=J_G^{-1}
{\iota^*}^{-1}(0)$ has dimension at most $p+n$. We will show
that there exists $z\in J_A^{-1}(0)$ with all $z_i\neq 0$.
If so, then $y=J_G(z)\in L^\perp$ and $J_G(z)$ is away
from the boundary of $\pi^*\Delta$ by Lemma \MomentImage.
Hence as $y$ vary slightly but arbitrarily
in $L^\perp$ we can always find $z$ with all $z_i\neq 0$
such that $y=J_G(z)$.
Moreover, $G=(S^1)^p$ obviously acts freely on
a small neighborhood of $z$ with all $z_i\neq 0$.
This shows that the dimension of $J_A^{-1}(0)$ at $y$
is at least $dim~G+dim~L^\perp=p+n$. This shows that $dim~J_A^{-1}(0)$
is exactly $p+n$. Since $G$ acts effectively on the zero set,
so does $A\subset G$. This also shows that the induced $T$ action
on $X_\Delta$ is effective, proving our third assertion.

We now show the existence of $z$. Pick $x$ in the interior of $\Delta$
so that $\bra x,u_i\ket>\lambda_i$, $\forall i$; hence
$y=(\bra x,u_1\ket,..,\bra x,u_p\ket)\in{\R^p}^*_>+\lambda$
is in the interior of $\pi^*\Delta$. Note that
$y\in L^\perp={\iota^*}^{-1}(0)$.
Pick any $z\in\C^p$
such that $0<y_i-\lambda_i=\half|z_i|^2$ for all $i$.
Then $J_G(z)=y\in{\iota^*}^{-1}(0)$, hence $z\in J_A^{-1}(0)$. $\Box$

\lemma{
Consider the moment map
$J_T:X_\Delta\ra\gt^*={\R^n}^*$, $[z]\mapsto (\pi^*)^{-1}J_G(z)$.
We have $J_T(X_\Delta)=\Delta$.}
\proof
Note that $J_T$ is well-defined because $J_G$ is $G$-equivariant
and that $z\in J_A^{-1}(0)\LRa J_G(z)\in Im~\pi^*$
by Lemma \ZeroSet, and that $\pi^*$ is injective. We have
$$
\pi^*J_T(X_\Delta)=J_G(J_A^{-1}(0))
=J_G(J_G^{-1}(\pi^*\Delta))
\subset\pi^*\Delta.
$$
Conversely, given $y\in\pi^*\Delta=({\R^p}^*_\geq+\lambda)\cap
{\iota^*}^{-1}(0)$, we have $y=J_G(z)$ for some $z\in\C^p$.
It follows that $\iota^*\circ J_G(z)=0$, hence $z\in J_A^{-1}(0)$,
and $\pi^*J_T([z])=J_G(z)=y\in\pi^*\Delta$. This shows that
$\pi^*J_T(X_\Delta)\supset\pi^*\Delta$. This proves that
$J_T(X_\Delta)=\Delta$. $\Box$

\theorem{For any $n$ dimensional rational polytope $\Delta\subset{\R^n}^*$,
the space $X_\Delta$ is a reduced Delzant $T$-space equipped
with the moment map $J_T$ and whose moment polytope
is $\Delta$.}
\proof
This follows from  the three preceding lemmas.
$\Box$

We now show that if $\Gamma\subset T$ is any finite subgroup, then $X_\Delta/\Gamma$
is a reduced Delzant $T/\Gamma$-space with the same moment polytope $\Delta$
in $Lie(T/\Gamma)^*=Lie(T)^*={\R^n}^*$. Be warned, however, that $X_\Delta/\Gamma$ is not homeomorphic to $X_\Delta$ in general. Note that we can also regard $X_\Delta/\Gamma$ as a $T$-space, but one in which $T$ does not act effectively.
We now prove the following more general assertion.

\theorem{If $X$ is a reduced $T$-space and $\Gamma\subset T$ is a finite subgroup, then
$X/\Gamma$ is a reduced $T/\Gamma$-space with the same moment polytope in $Lie(T/\Gamma)^*=Lie(T)^*$.}
\proof
Let $1\ra A\ra G\ra T\ra 1$ be an exact sequence and
suppose that $X$ is obtained by reducing the $G$-manifold $(M,\omega)$  with respect to $A$.
We have another exact sequence $1\ra A'\ra G\ra T/\Gamma\ra1$
where $A'$ is the kernel of the composition map $G\ra T\ra T/\Gamma$.
We can reduce $(M,\omega)$ with respect to $A'$. Call the reduced space $X'$. Since $A$ and $A'$ have the same moment map, and $A'$ is the preimage of $\Gamma$ under $G\ra T$, we have an exact sequence $1\ra A\ra A'\ra\Gamma\ra1$. Thus we can apply reduction in two stages
\SL~ and obtain $X'=X/\Gamma$, which is a $T/\Gamma$-space with moment map
$J_{T/\Gamma}=J_T$. Since $J_T$ is $T$-invariant,
it follows that $J_T(X/\Gamma)=J_T(X)$. $\Box$

In particular, if we let $\Gamma=\Z^n/\Pi$, then the preceding two theorems imply that there exists a reduced Delzant $\R^n/\Z^n$-space ($T/\Gamma=\R^n/\Z^n$) whose moment polytope is $\Delta$.

\subsec{From reduced spaces to complete toric varieties}

We continue to use the notations of the preceding section.

Let $F$ be the relative interior of a face of $\Delta$. Then there is a unique
index set $I_F\subset\{1,..,p\}$ such that the closure
$$
\bar F=\{x\in\Delta|\bra u_i,x\ket=\lambda_i,~i\in I_F\}.
$$
Put
$$
V_F:=\{z\in\C^p|z_i=0\LRa i\in I_F\}\subset\C^p
$$
which has complex dimension $p-\#I_F$ and
is clearly invariant under the action of $G_\C=(\C^\times)^p$.
Since the $I_F$ are pairwise distinct, the $V_F$
are pairwise disjoint, and
$$
V:=\cup_F V_F\subset\C^p
$$
defines a $G_\C$-invariant disjoint union of $G_\C$-orbits. It is not hard to show that this is also a complex stratified space. Moreover this stratification corresponds to the stratification of $\Delta$ by the interiors of its faces $F$. A note about notation: since we will not be discussing the local normal form in this subsection, $V$ here should not be confused with the symplectic slice in $M$.

\lemma{$J_G^{-1}(\pi^*F)\subset V_F$.
In particular $J_A^{-1}(0)\subset V$. Moreover
the inclusion is $G$-equivariant.}
\proof
The last assertion is obvious.
The second assertion follows from that the first assertion and that
$J_A^{-1}(0)=J_G^{-1}(\pi^*\Delta)=\cup_F J_G^{-1}(\pi^* F)$.
For the first assertion, note that
$J_A(z)=0\LRa J_G(z)\in\pi^*F,~\exists!F$.
This holds iff $J_G(z)=\pi^*(x),~\exists~x\in F$.
But $x\in F$ means that $\bra x,\pi^*(e_i)\ket=\lambda_i
\LRa i\in I_F$. It follows that
$\bra J_G(z),e_i\ket=\bra\pi(x),e_i\ket=\lambda_i\LRa i\in I_F$.
Since $J_G(z)=\half(|z_1|^2,..,|z_p|^2)+\lambda$,
this means that $z_i=0\LRa i\in I_F$, i.e. $z\in V_F$. $\Box$

\corollary{Consider $J_T:J_A^{-1}(0)\ra{\R^n}^*=\gt^*$. Then
$J_T^{-1}(F)=J_G^{-1}(\pi^*F)=V_F\cap J_A^{-1}(0)$.}
\thmlab\OverEachFace
\proof
The second equality follows from the preceding lemma
and the fact that $J_A^{-1}(0)$ is the disjoint
union of the $J_G^{-1}(\pi^*F)$ as $\bar F$ ranges over
faces of $\Delta$. For the first equality, consider
the commutative diagram
$$\matrix{
&J_A^{-1}(0) &\subset & \C^p& &\cr
&J_T\da & & \da J_G& &\cr
\Delta\subset&{\R^n}^*&\br\pi^*\over\hra &{\R^p}^*&\br\iota^*\over\ra& L_\R^*.}
$$
For $z\in\C^p$, we have
$$
J_T(z)\in F\LRa\pi^* J_T(z)=J_G(z)\in\pi^*F\LRa
z\in J_G^{-1}(\pi^*F).
$$
This completes the proof. $\Box$

Recall that $\sigma_F:=\sum_{i\in I_F}\R_\geq u_i\subset\R^n$,
so that $\sigma_F^\vee=\{x\in{\R^n}^*|i\in I_F\Lra\bra x,u_i\ket\geq0\}$.
It follows that all $z_j^{\bra x,u_j\ket}$ are well-defined
for any $z\in V_F$ and $x\in\sigma_F^\vee\cap\Pi^*$.
For $F\neq\Delta$, we define
$$
\phi_F:V_F\ra U_F:=Hom_{sg}(\sigma_F^\vee\cap\Pi^*,\C),
~~~z\mapsto (x\ra\prod_j z_j^{\bra x,u_j\ket}).
$$
Here $\C$ is regarded as a multiplicative semigroup.
For $F=\Delta$, we use the same definition except that $\C$
is replaced by $\C^\times$. Let $\Sigma=\Sigma_\Delta$ be the fan consisting of the cones $\sigma_F$, as $\bar F$ ranges over all faces of $\Delta$.
By definition, the toric variety $\P_\Sigma(\Pi)$ is the union of
the $U_F$ modulo the relations given by the open embeddings
$U_F\hra U_E$, $\bar E\subset \bar F$.
Gluing the maps $\phi_F$ together, we get a map
$$
\phi:V\ra\P_\Sigma(\Pi)=\cup_F U_F
$$
where we identify $U_F$ with a subset of $\P_\Sigma(\Pi)$.
It is helpful to keep in mind that
$$
\bar E\subset \bar F\LRa I_E\supset I_F\LRa\sigma_E\supset\sigma_F
\LRa\sigma_E^\vee\subset\sigma_F^\vee
\LRa U_E\supset U_F.
$$

The group $U_\Delta:=Hom(\Pi^*,\C^\times)$ acts on each $U_F$,
by multiplication of functions:
$U_\Delta\times U_F\ra U_F,~f,g\mapsto f\cdot g$.
In turn, the group $T_\C:=\C^n/\Pi$ acts on
$U_F$ via the group isomorphism
$T_\C\ra U_\Delta$, $u\mapsto e^{2\pi i\bra-,u\ket}$.
The group $G_\C=\C^p/\Z^p$ acts on $U_F$ via
the group homomorphism $\pi:G_\C\ra T_\C$
induced by $\pi:\Z^p\ra\Pi$, $e_i\mapsto u_i$.
Explicitly, this action $G_\C\times U_F\ra U_F$
is given by $(g,v+\Z^p)\mapsto e^{2\pi i\bra -,\pi(v)\ket}\cdot g$.
Note that $A_\C:=L_\C/L=Ker(\pi:G_\C\ra T_\C)$.

\lemma{$\phi:V\ra\P_\Sigma(\Pi)$ is continuous, $G_\C$-equivariant and $A_\C$-invariant.}
\proof
The argument for continuity is standard.
To show the $G_\C$-equivariance, note that $G_\C$ acts on $\C^p$ hence on $V$, by
$f\cdot z=(e^{2\pi i\bra e_1^*,v\ket} z_1,..,e^{2\pi i\bra e_p^*,v\ket} z_p)$
for $f=v+\Z^p$.
It follows that for $z\in V_F$,
$$\eqalign{
\phi_F(g\cdot z):x&\mapsto
\prod_j(e^{2\pi i\bra e_j^*,v\ket}z_j)^{\bra x,u_j\ket}\cr
&=\prod_j(e^{2\pi i\bra e_j^*,v\ket}z_j)^{\bra x,\pi(e_j)\ket}\cr
&=\prod_j(e^{2\pi i\bra e_j^*,v\ket}z_j)^{\bra \pi^*(x),e_j\ket}\cr
&=e^{2\pi i\bra \pi^*(x),v\ket}\prod_jz_j^{\bra x,u_j\ket}\cr
&=e^{2\pi i\bra x,\pi(v)\ket}\prod_jz_j^{\bra x,u_j\ket}\cr
&=g\cdot\phi_F(z)(x).
}$$
This shows that $\phi_F$ is $G_\C$-equivariant.
Finally since $\pi(L)=0$, the last formula shows that that
for $g\in A_\C=L_C/L$ we have $\phi_F(g\cdot z)=\phi_F(z)$,
hence $\phi_F$ is $A_\C$-invariant.
$\Box$

\corollary{$\phi$ induces a $T$-equivariant continuous map $\phi:X_\Delta=J_A^{-1}(0)/A\ra
\P_\Sigma(\Pi)$.}

\lemma{(Transversality) Consider the vector part $L_\R$
of the noncompact group $A_\C=L_\C/L=L_\R\times A$
acting on $\C^p$. Each $L_\R$-orbit is transversal to the zero set $J_A^{-1}(0)$.
In other words, for $v\in L_\R$ and $z\in J_A^{-1}(0)$
we have $v\cdot z\in J_A^{-1}(0)\Lra v\cdot z=z$.}
\proof
The $L_\R$-action on $\C^p$ is
$z\mapsto v\cdot z=(e^{v_1}z_1,..,e^{v_p}z_p)$ for $v\in L_\R\subset\R^p$
i.e. $\pi(v)=0$. Since $v\in L_\R$ and $J_A:\C^p\ra L_\R^*$,
the following function $g:\R_>\ra\R$ makes sense:
$$\eqalign{
g(s)&=\bra J_A(s^v\cdot z),v\ket\cr
&=\sum_i(\half|z_i|^2s^{2v_i}+\lambda_i)\bra\iota^*(e_i^*),v\ket\cr
&=\sum_i(\half|z_i|^2s^{2v_i}+\lambda_i)v_i.
}$$
Since $z,v\cdot z\in J_A^{-1}(0)$, it follows that
that $g(1)=g(e)=0$. Thus $g'(s)=0$ for some $s\in[1,e]$.
Hence
$$
g'(s)={1\over s}\sum_i|z_i|^2 s^{2v_i}v_i^2=0
$$
implying that $|z_i|v_i=0$ for all $i$. This shows that
$e^{v_i}z_i=z_i$ for all $i$, i.e. $v\cdot z=z$. $\Box$

\corollary{Let $z\in J_A^{-1}(0)$.
Then $(A_\C\cdot z)\cap J_A^{-1}(0)=A\cdot z$.
In particular the inclusion $J_A^{-1}(0)\subset V$ induces
a $T$-equivariant inclusion map $\tilde\psi:J_A^{-1}(0)/A\hra V/A_\C$.}

\lemma{(Open Embedding) Put $X^\circ:=J_G^{-1}(\pi^*\Delta^\circ)$.
The map $\psi:X^\circ\times L_\R\ra V_{\Delta^\circ}$, $(z,v)\mapsto v\cdot z$,
is an $A_\C$-equivariant open mapping. In particular,
$\psi(X^\circ/A)\subset V_\Delta/A_\C$ is open.}
\proof
Given the first assertion, $(X^\circ\times L_\R)/A_\C
=X^\circ/A$ is mapped to an open set in $V_\Delta/A_\C$
yielding the second assertion.
Note that $A_\C=A\times L_\R$ acts on $X^\circ\times L_\R$ by
$(t,v),(z,v')\mapsto(t\cdot z,v+v')$, and on $V_\Delta\subset\C^p$ by
$(t,v),z\mapsto tv\cdot z$. Hence $\psi$ is $A_\C$-equivariant.
To prove the openness of $\psi$, it suffices to show
that $X^\circ\subset V_\Delta$ is an embedded closed submanifold
and each $L_\R$-orbit in $V_\Delta$ that meets $X^\circ$
meets it transversally, which follows from the Transversality Lemma
because $X^\circ\subset J_A^{-1}(0)$.
Now by Corollary \OverEachFace,
$$
X^\circ=V_\Delta\cap J_A^{-1}(0)
=\{z\in{\C^\times}^p|J_A(z)=0\}.
$$
We have $J_A(z)=0\LRa\half|z|^2+\lambda\in L_\R^\perp$ in $\R^p$.
Thus $X^\circ$ is defined by a finite set of real quadratic equations
in ${\C^\times}^p$. If we choose a basis
$l_1,..,l_k$ of $L_\R$, then the equations are
$(\half|z|^2+\lambda)\cdot l_j=0$. This shows that
$X^\circ\subset V_\Delta$ is closed. The normal vectors
of those hypersurfaces are
$$
(\bar z_1l_{1j},..,\bar z_pl_{pj},
z_1l_{1j},..,z_pl_{pj}),~~~j=1,..,k.
$$
Putting them in a $2p\times k$ matrix, we get
$\left[\matrix{diag(\bar z)l\cr diag(z)l}\right]$ where $l=(l_{ij})$.
Since $z\in(\C^\times)^p$, this matrix and $l$ have
the same rank i.e. $k$. This shows that the complete
intersection $X^\circ\subset V_\Delta$ of those real hypersurfaces is
an embedded smooth submanifold of dimension $2p-k$. $\Box$

For $v\in F$,  we have $\R(F-v)=\cap_{i\in I_F} u_i^\perp=\sigma_F^\perp$,
so that $dim~\sigma_F^\perp=dim~F$.
Since $\sigma_F^\perp\cap\Pi^*$ is primitive sublattice of $\Pi^*$
of rank $dim~\sigma_F^\perp=dim~F$, it follows that
$$
W_F:=Hom(\sigma_F^\perp\cap\Pi^*,\C^\times)\cong(\C^\times)^{dim~F}.
$$
Note that $T_\C\cong Hom(\Pi^*,\C^\times)$ acts on $W_F$ by function multiplications,
and the action is clearly transitive.

\lemma{There is a $T_\C$-equivariant inclusion $W_F\hra U_F$
given by extension by zero. Thus we can regard $W_F\subset U_F$.}
\proof
Note that $\sigma_F^\perp$ is the largest vector subspace
in the cone $\sigma_F^\vee$ in ${\R^n}^*$. Let $\tau$ be any
finitely generated cone over $\R$, and $\tilde\tau$
the largest vector subspace in $\tau$.
Define a projection map $f:\tau\ra\tilde\tau$,
by $f(x)=x$ if $x\in\tilde\tau$, and $f(x)=0$ if $x\notin\tilde\tau$.
Using the fact that
$\forall x,y\in\tau,~x+y\in\tilde\tau\LRa x,y\in\tilde\tau$,
it's straightforward to check that $f$
is a real semi-group homomorphism. We apply this
to the case $\tau=\sigma_F^\vee$ and $\tilde\tau=\sigma_F^\perp$. Define
$W_F\ra U_F=Hom(\sigma_F^\perp\cap\Pi^*,\C^\times)$, $x\mapsto x\circ f\in U_F$. It is an inclusion because $\sigma_F^\perp\subset\sigma_F^\vee$,
and so $x\circ f=x'\circ f$ implies that $x=x'$ on $\sigma_F^\perp$.
The map is clearly $T_\C$-equivariant. 
$\Box$

\lemma{$\phi_F(V_F)\subset W_F$.}
\proof
Recall that $\phi_F(V_F)\subset U_F$ and that for
$\alpha\in U_F$, we have
$$
\alpha\in W_F\LRa\alpha(x)=0,~\forall x\notin
\sigma_F^\perp\cap\Pi^*.
$$
Let $z\in V_F$ so that $z_i=0\LRa i\in I_F$.
Then $\phi_F(z)(x)=\prod_i z_i^{\bra x,u_i\ket}=0$
iff $\bra x,u_i\ket>0$, $\exists i\in I_F$,
iff $x\notin\sigma_F^\perp\cap\Pi^*$. If follows
that $\phi_F(z)\in W_F$. $\Box$

\lemma{By dropping all the zeros in $V_F$, view
$V_F\cong(\C^\times)^{p-\#I_F}$ as a group. Then the map
$\phi_F:V_F\ra W_F$ is an $A_\C$-invariant group homomorphism.
Moreover the $A_\C$-action on $V_F$ is equivalent
to the one obtained from the
group homomorphism $A_\C=L_\C/L\ra V_F$, $v+L\mapsto
e^{2\pi iv}\cdot\One_F$ where $\One_F\in V_F$ is the
unique vector with entries 0 or 1.}
\proof
Straightforward. $\Box$

\lemma{For each $\alpha\in W_F$, the fiber $\phi_F^{-1}(\alpha)$
in $V_F$ is a single $A_\C$-orbit.}
\proof
By the preceding lemma, it suffices to show that
$\phi_F^{-1}(1)=A_\C\cdot\One_F$. The inclusion $\supset$
is obvious from the definition of $\phi_F$.
Suppose $z\in\phi_F^{-1}(1)$ in $V_F$.
This means that $z_i=0\LRa i\in I_F$ and that
$$
\prod_j z_j^{\bra x,u_j\ket}=1,~~~
\forall x\in R:=\sigma_F^\perp\cap\Pi^*.
$$
We want to find $v\in L_\C$ such that $e^{2\pi iv}\cdot\One_F=z$,
i.e. $z_j=e^{2\pi i v_j}$ for $j\notin I_F$. Pick any $w\in\C^p$ such that
$0\neq z_j=e^{2\pi i w_j}$ for $j\notin I_F$.
Note that we are free to change
the values of the $w_j$, $j\in I_F$, without changing this relation.
We have
$$
1=\prod_j e^{2\pi i\bra x,w_j u_j\ket},~~~~x\in R.
$$
Since $\bra x,u_j\ket=0$, $\forall j\in I_F$, this
holds for any values we choose for $w_j$, $j\in I_F$.
This equation now reads $1=e^{2\pi i\bra x,\pi(w)\ket}$
where $w=\sum_j w_j e_j\in\C^p$, because $\pi(e_j)=u_j$.
It follows that $\pi(w)\in R^*=(\sigma_F^\perp\cap\Pi^*)^*
=\Pi/\Pi_F$, where $\Pi_F:=\R\sigma_F\cap\Pi$. This means that
there exists $u\in\R\sigma_F=\sum_{i\in I_F}\R u_i
=\pi(\sum_{i\in I_F}\R e_i)$.
such that $\pi(w)+u\in\Pi$. In other words, we can change
the values of the $w_j$, $j\in I_F$, so that $\pi(w)\in\Pi=\pi(\Z^p)$.
Finally, pick a vector $w'\in\Z^p$ so that $\pi(w+w')=0$.
Then $v=w+w'\in L_\C$, and we have $z_j=e^{2\pi iw_j}=e^{2\pi iv_j}$
for $j\notin I_F$ because $w'\in\Z^p$. This completes the proof. $\Box$

\lemma{The map
$\phi:V\ra\P_\Sigma(\Pi)$ induces a continuous $T$-equivariant bijection
$\tilde\phi:V/A_\C\ra\P_\Sigma(\Pi)$.}
\proof
This follows from the preceding lemma
and the fact that the $V_F$ are
pairwise disjoint in $V$ and that the $W_F$ are pairwise
disjoint in $\P_\Sigma(\Pi)$. $\Box$

\theorem{The natural maps
$X_\Delta{\br\tilde\psi\over\ra} V/A_\C
{\br\tilde\phi\over\ra}\P_\Sigma(\Pi)$ are $T$-equivariant homeomorphisms.}
\proof
The maps are $T$-equivariant by constructions.
Since $\P_\Sigma(\Pi)$ is Hausdorff (in
the usual analytic topology), the
preceding lemma implies that $V/A_\C$ is Hausdorff.
Since $X=X_\Delta=J_A^{-1}(0)/A$ is compact and $\tilde\psi$
is continuous, it follows that $\tilde\psi(X)\subset V/A_\C$
is compact, hence closed.
In particular $\tilde\psi(X)\cap V_{\Delta^\circ}/A_\C$ is closed
in $V_{\Delta^\circ}/A_\C$.

Since the $V_F/A_\C$ are pairwise disjoint,
by Corollary \OverEachFace~ we have
$\tilde\psi(X^\circ)=\tilde\psi(J_G^{-1}(\pi^*\Delta^\circ))
=\tilde\psi(X)\cap V_{\Delta^\circ}/A_\C$.
By the Open Embedding Lemma, this is open in $V_{\Delta^\circ}/A_\C$.
But $V_{\Delta^\circ}/A_\C$ is an $n$ dimensional algebraic torus, hence connected.
It follows that $\tilde\psi(X^\circ)$ must be all of $V_{\Delta^\circ}/A_\C$.
In particular it is dense in $V/A_\C$. Since $\tilde\psi(X)$
is also closed in $V/A_\C$ and contains
the dense subset $\tilde\psi(X^\circ)$, it must be all of $V/A_\C$.
By the corollary to the Transversality Lemma,
it follows that $\tilde\psi$ is a continuous bijection.
Since $\tilde\phi$ is also a continuous bijection, so is the composition $\tilde\phi\circ\tilde\psi$.
Since both $X$ and $\P_\Sigma(\Pi)$ are finite union of manifolds,
this composition is a homeomorphism. Thus both $\tilde\psi,\tilde\phi$ are homeomorphisms.
$\Box$

\newsec{Classification of Reduced Delzant Spaces}

In this section, we reconstruct the $T$-equivariance homeomorphism type of a reduced Delzant $T$-space $X$ from its moment polytope in $\gt^*$. Again, the notations introduced in section 2 shall remain in force here.

Fix a rational polytope $\Delta$ in $\gt^*$, and let $F$ be the (relative) interior of a given face of $\Delta$. Put $X_F=J_T^{-1}(F)\subset X$, $Z_F=J_G^{-1}(0\times F)\subset M$.
By a corollary of the Face Theorem, the points in $Z_F$ have the
same stabilizer type $H\subset G$. Since $F$ is smoothly contractible, by choosing a base point in $F$ we have a $G$-equivariant isomorphism $Z_F\cong G/H\times F$ over $F$.
Taking the $A$-orbit spaces, we get a $T$-equivariant isomorphism
$X_F\cong T/B\times F$. We shall make the identification of spaces by means of these isomorphisms.

For each $F$, we would like to describe a neighborhood of $X_F$ in $X$,
as a stratified symplectic space, and the moment map $J_T$ there in order
to reconstruct the topology of $X$ and its moment map. The idea is to
first describe a neighborhood of $Z_F$ in $M$ by using the constant rank embedding theorem of Sjamaar-Lerman and the minimal coupling procedure of Sternberg and Weinstein. Then we reduce that neighborhood with respect to $A$. From this we get a kind of {\it semi-global form} in a neighborhood of each stratum of the reduced space. Some of the machinery used in this section is borrowed from \SL.
Finally we reconstruct $X$ and its moment map $J_T$ by gluing together these semi-global forms.

\subsec{Geometric normal bundle}

In this subsection, we shall write $Z=Z_F$.
By the corollary to Lemma \StabilizerType, we saw that the two-form $\tau:=\omega|Z$ has constant rank. By the constant rank embedding theorem, $Z\hra M$
corresponds to a symplectic vector bundle $\cN\ra Z$. In fact
$\cN$ is $G$-equivariantly isomorphic to $G\times_H V\times F$,
where $V$ is the symplectic slice at a chosen point $\tilde p\in G/H\times F$. Note that $Z$ sits inside $G\times_HV\times F$ as the set $(G\times_H0)\times F$. Note that by construction of the local normal form, we can assume that $V$ depends only on $v\in F$.

Let $\cV$ the subbundle of $TZ$ whose typical fiber is the radical of the two-form $\tau$, i.e.
$$
\cV_z=\{v\in T_zZ|\tau(v,u)=0,~\forall u\in T_zZ\}.
$$ 
Let $\pi_{\cV^*}:\cV^*\ra Z$ be the dual bundle. Then the geometric normal bundle of $Z\hra M$
is isomorphic to the Whitney sum $\cV^*\oplus\cN$.
Since $\omega$ is a $G$-invariant form, $\cV$ is also a $G$-equivariant bundle.

\lemma{For $(gH,v)\in Z$, we have $\cV_{(gH,v)}=T_{gH} (gD/H)\times 0$ where
the right hand side is viewed as the fiber
of a subbundle of $T(G/H\times F)=T(G/H)\times(F\times\gn^*)$ over $Z$.}
\thmlab\VLemma
\proof
In the local normal form of $Z$, the symplectic form $\omega$ is given by
the form inherited from a neighborhood of the zero section of $G\times_H(\gq^*\times V) \subset(G\times\gg^*\times V)/H$, where the symplectic form is the form on $G\times\gg^*$
plus $\omega_V$ (cf. Remark \SSOnLocalNormalForm.) The constant rank form $\tau$ is the restriction to the submanifold $Z\equiv G/H\times F\subset G/H\times(pt+\gn^*)$, where $pt$ is any point in $F$.
The tangent space at $(gH,v)$ of $Z$ is $T_{gH}(G/H)\times\gn^*\cong\gg/\gh\times\gn^*$.
The bilinear form $\tau$ on this vector space is inherited from the pairing on $\gg\times\gg^*$.
Thus the radical of $\tau$ on $T_{(gH,v)}Z$ is $T_{gH}(gD/H)\cong \gd/\gh$. $\Box$

Choose a $G$-equivariant splitting $TZ=\cV\oplus\cH$.
This corresponds to a choice of $G$-equivariant section $s:\cV^*\ra T^*Z$ of the bundle map
 $T^*Z\twoheadrightarrow\cV^*$ over $Z$.
Put
\eqn\MuDefn{
\mu=(\pi_{\cV^*})^*\tau+s^*\gamma
}
where $\gamma$ is the canonical symplectic form on $T^*Z$. The two-form $\mu$ on $\cV^*$
is $G$-invariant, closed, and non-degenerate near the zero section
$Z\hra\cV^*$. 

Fix $t\in T$ and $g\in\pi^{-1}(t)$.

\lemma{The fiber of the composed map
$\cV^*\ra Z\ra X_F$ at $(tB,v)\in X_F$ is the manifold $T^*(gD/H)\times v$.
Moreover the $\mu$ restricted to this fiber
agrees with the canonical form on $T^*(gD/H)$.}
\thmlab\RestrictionofMu
\proof
The fiber of $Z=G/H\times F\ra X_F=T/B\times F$ at $(tB,v)$
is $gD/H\times v$ where $\pi:g\mapsto t$ under the map $\pi:G\ra T$. It follows that the fiber of the composed map is the bundle $\cV^*|(gD/H\times v)$. By Lemma \VLemma,
the typical fiber of this bundle at $(gdH,v)$, where $d\in D$, is $T^*_{gdH}(gdD/H)\times v=T^*_{gdH}(gD/H)\times v$. 
It follows that the fiber over $(tB,v)$ of the bundle $\cV^*\ra X_F$ is $T^*(gD/H)\times v$.
The typical tangent space of this fiber is $T_{gdH}(gD/H)\times T_{gdH}^*(gD/H)$.

First we claim that $(\pi_{\cV^*})^*\tau$ restricted to this fiber is zero, i.e.
that as a bilinear form on the tangent space
$T_{gdH}(gD/H)\times T_{gdH}^*(gD/H)$, it is identically zero. 
Here the first factor is $\cV_{(gdH,v)}$, which is in the
radical of $\tau$ in $T_{(gdH,v)}Z$ hence $(\pi_{\cV^*})^*\tau$
evaluated on this factor is zero. The second factor is
$\cV_{(gdH,v)}^*\hra T_{(gdH,v)}^*Z$ is the fiber of the vector bundle $\cV^*$. But $(\pi_{\cV^*})^*\tau(\cV_{(gdH,v)}^*,-)=0$. This proves our claim.
 
 It remains to show that $s^*\gamma$ as a bilinear form on
 $T_{gdH}(gD/H)\times T_{gdH}^*(gD/H)$ agrees with the canonical form.
By definition $\gamma$ as a bilinear form on a typical tangent space of $T^*Z$
 is given by the pairing on $T_zZ\times T_z^*Z$, $z=(gdH,v)$.
 We have $TZ=\cV\oplus\cH$ and $T^*Z=\cV^*\oplus\cH^*$ via the section $s$.
 Thus the pairing on $T_zZ\times T^*_zZ$ just restricts to the pairing on
 $\cV_z\times\cV_z^*$.  $\Box$
 
 \lemma{The fiber of the composed map
$\cN\ra Z\ra X_F$ at $(tB,v)\in X_F$ is the manifold 
$gD\times_HV\times v\subset G\times_HV\times F=\cN$.}
\thmlab\NOverXF
\proof
Again, the fiber at $(tB,v)$ of $Z\ra X_F$ is $gD/H\times v\subset G/H\times F$.
Since $\cN\ra Z$ is the bundle $G\times_HV\times F$, its restriction to $gD/H\times v$
is clearly given by our assertion. $\Box$

\subsec{Minimal coupling procedure}

In this subsection, we construct the symplectic structure near the zero section of the geometric normal bundle $\cN\oplus\cV^*$ of $Z$ in $M$.
Consider the symplectic normal bundle
$\cN=G\times_H V\times F$ over $Z$ with typical fiber $V$. We can view the symplectic
structure $\omega_V$ as the imaginary part of an $H$-invariant hermitian form on $V$,
and make $\cN$ a hermitian vector bundle. Then the symplectic $H$-action now becomes a unitary action. Put $U=U(V)$, the unitary group on $V$, and let 
$\tilde\pi_\cN:Fr(\cN)\ra Z$ be the unitary frame bundle of $\cN$,
which is a principal $U$-bundle on which $U$ acts on the right. The left $G$-action on $\cN$ makes $Fr(\cN)$ an $G$-equivariant bundle.

Let $\pi_{\cN^\#}:\cN^\#\ra\cV^*$ be the pullback of the bundle $\cN\ra Z$ along $\cV^*\ra Z$.
Then we have the $G$-equivariant commutative diagram:
\eqn\NDiagram{\matrix{
\cN^\# & \longrightarrow &\cN\cr
\pi_{\cN^\#} \da\hskip.3in & &\hskip.3in \da \pi_{\cN}\cr
\cV^* &\br\pi_{\cV^*}\over\longrightarrow & Z.
}}
Note that $\cN^\#=\cN\oplus\cV^*$ as bundles over $Z$.
The pullback of $Fr(\cN)$ along $\cV^*\ra Z$ is the principal $U$-bundle 
$Fr(\cN^\#)$ over  $\cV^*$. We have another $G$-equivariant commutative diagram
\eqn\FrameDiagram{\matrix{
Fr(\cN^\#) &\br\tilde\pi_{\cV^*}\over\longrightarrow &Fr(\cN)\cr
\tilde\pi_{\cN^\#} \da\hskip.3in & &\hskip.3in \da \tilde\pi_{\cN}\cr
\cV^* &\br\pi_{\cV^*}\over\longrightarrow & Z.
}}
Fix a $G\times U$-invariant connection $\theta$ on $Fr(\cN)$, and let
$\tilde\theta=(\tilde\pi_{\cV^*})^*\theta$. Then there exists a $G\times U$-invariant symplectic form
(see section 8 \SL)
$$
\sigma={(\tilde\pi_{\cN^\#}})^*\mu-d\bra pr_2,\tilde\theta\ket
$$
in a neighborhood of $Fr(\cN^\#)\times 0$ in $Fr(\cN^\#)\times\gu^*$.
Here $pr_2:Fr(\cN^\#)\times\gu^*\ra\gu^*$ is the projection and is also the negative of a moment map for the Hamiltonian $U$-action; $\mu$ is defined in \MuDefn.

The $U$-action on a neighborhood in $Fr(\cN^\#)\times\gu^*\times V$
is Hamiltonian with a moment map $(l,\eta,\nu)\mapsto\Psi(\nu)-\eta\in\gu^*$,
where $\Psi$ is a $U$-moment map on $V$. Let $(Fr(\cN^\#)\times\gu^*\times V)_0$
be the zero set of the moment map. Then we have a natural $G$-equivariant isomorphism
$$
\cN^\#=Fr(\cN^\#)\times_UV{\br j\over\ra}(Fr(\cN^\#)\times\gu^*\times V)_0/U.
$$
Define
\eqn\sigmaSharp{
\sigma^\#=j^*(\sigma+\omega_V)'
}
where $(\sigma+\omega_V)'$ is the form induced on the reduced space $(Fr(\cN^\#)\times\gu^*\times V)_0/U$. Then $\sigma^\#$ is a $G$-invariant symplectic form in a neighborhood of
the zero section $\cV^*\hra\cN^\#$.

This symplectic structure on $\cN^\#$ as a bundle
over $X_F=T/B\times F$ can also be easily described fiberwise. Recall that the fiber at $(tB,v)\in X_F$ of the composed map $\cV^*\ra Z\ra X_F$ is $T^*(gD/H)\times v=gD\times_H(\gd/\gh)^*\times v$ 
and that its canonical form coincides
with $\mu$ restricted to this fiber by Lemma \RestrictionofMu. Here $t=\pi(g)$ under $\pi:G\ra T$.
Also the fiber of the map $\cN\ra X_F$ is $gD\times_HV\times v$.
Thus the fiber at $(tB,v)$ of the bundle $\cN^\#=\cV^*\oplus\cN\ra X_F$ is
\eqn\FiberOverXF{
gD\times_H((\gd/\gh)^*\times V)\times v.
}
Note that this is canonically isomorphic to the reduced space  $(T^*gD\times V)_0/H$
where $T^*gD$ is given the canonical form and $V$ is given $\omega_V$.
By direct calculations, we find 

\lemma{The induced symplectic structure on the reduced space $(T^*gD\times V)_0/H$
coincides with $\sigma^\#$ restricted to the fiber \FiberOverXF~ of $\cN^\#\ra X_F$.}

Here is a schematic picture of the structures we have gathered so far:
$$\matrix{
Fr(\cN^\#)\times\gu^*,\sigma & \longrightarrow & Fr(\cN^\#),\tilde\theta &
\br\tilde\pi_{\cV^*}\over\longrightarrow & Fr(\cN),\theta & &\cr
pr_2\da& &\tilde\pi_{\cN^\#}\da\hskip.3in & &\hskip.3in \da\tilde\pi_{\cN} & &\cr
\hskip.2in\gu^*& &\cV^*,\mu&\br\pi_{\cV^*}\over\longrightarrow & Z ,\tau & \hra & M,\omega\cr
& &s\da\hskip.1in & & \da &  &  \cr
& & T^*Z,\gamma & & X_F & & \cr
}$$
$$\eqalign{
\mu&=(\pi_{\cV^*})^*\tau+s^*\gamma\cr
\sigma&=(\tilde\pi_{\cN^\#})^*\mu-d\bra pr_2,\tilde\theta\ket\cr
\sigma^\#&=j^*(\sigma+\omega_V)'\cr
\tilde\theta&=(\tilde\pi_{\cV^*})^*\theta\cr
\cN^\#&=Fr(\cN^\#)\times_U V{\br j\over\ra}(Fr(\cN^\#)\times\gu^*\times V)_0/U\cr
& ~~~~[e,\nu]\mapsto[e,\Psi(\nu),\nu].
}$$

\subsec{$G$-moment map on $Z_F$}

We want to analyze the $G$-moment map $J_G:M\ra\gg^*$ modeled on a neighborhood of
the zero section of $Z\hra\cN^\#$. This is where we get the crucial local information we need about our $G$-moment map. We shall write $\cN^\#$ when we
mean a neighborhood of its zero section, and denotes the moment map here by $J^\#$.

Fix a point $(tB,v)\in X_F$ and consider the fiber \FiberOverXF~ of the bundle $\cN^\#\ra X_F$. The normal subgroup $D\subset G$ acts on the left. The $D$-moment map is given by
$$
Dg\times_H((\gd/\gh)^*\times V)\ra\gd^*,~~~ [dg,\eta,\nu]\mapsto Ad^*(dg)(\eta+\Phi(\nu))+c
$$
where $c\in\gd^*$ is a constant to be determined, and $\Phi$ is the $H$-moment map on $V$.
Now $\gg^*=\gn^*\oplus\gd^*$, as $Ad^*(G)$-module.  Hence $J^\#:\cN^\#\ra\gg^*$ has the shape $J^\#=(J_1,J_2)$ where $J_2$ is the moment map for the subgroup $D\subset G$
given above when restricted to the fiber \FiberOverXF.

\lemma{Let $\tilde q=(gH,v)\in Z\hra\cN^\#$. Write $v=(v_1,v_2)\in\gn^*\oplus\gd^*$. Then $J_i(\tilde q)=v_i$ for $i=1,2$. Moreover $J_2$ restricted to fiber \FiberOverXF~ is
 $$
J_2([dg,\eta,\nu]\times v)=Ad^*(dg)(\eta+\Phi(\nu))+v_2
$$}
\thmlab\RestrictionOfJ
\proof
The $G$-moment map $J^\#$ restricted to $Z=G/H\times F\subset M$ is the projection $pr_F:Z\ra F$. The point $\tilde q\in Z$ has the shape $[dg,0,0]\times v$ in \FiberOverXF. Since  $J^\#(\tilde q)=pr_F(\tilde q)=v$, the first assertion follows. 
On the other hand, the expression for $J_2$ above yields $c=J_2(\tilde q)=v_2$. $\Box$
 
We now analyze $J_1$. Recall that in Lemma \FiniteCoverN, we have a central torus
$N'$ in $G$ such that $N'\ra N$ is a finite cover under $G\ra N$. First we want to establish the existence of a horizontal distribution on the principal $U$-bundle $Fr(\cN^\#)$
which is $G\times U$-invariant and contains every tangent
space $T_x(N'\cdot x)$ of an $N'$-orbit. 
Since the bundles $\cN,\cV^*$ are all $G$-equivariant, it is enough to do this on 
$$
P=Fr(\cN)
$$
and then pullback the distribution to $Fr(\cN^\#)$. 

\lemma{$N'\times U$ acts locally freely on $P$.}
\proof
It is easy to check that $G$ acts freely on the principal $U$-bundle $P$.
Suppose $(n,u)\in N'\times U$ stabilizes $x\in P$,
i.e. $nxu^{-1}=x$, i.e. $nx=xu$. Project this down to the base $Z$; we get
$n\tilde\pi_\cN(x)=\tilde\pi_\cN(x)$ because $\tilde\pi_\cN$ is $U$-invariant and $G$-equivariant.
Hence $n\in Stab_{N'}\tilde\pi_\cN(x)$. But any stabilizer subgroup of $N'$
acting on $Z=G/H\times F$ is finite, because $\gn'\cap\gd=0$ and $\gh\subset\gd$
imply that $N'\cap H$ is finite. This shows that for a given $x$, $n$ ranges over
only a finite subset of $N'$. Now for each $n$ in that finite subset,
there is at most one $u\in U$ such that $(n,u)$ stabilizes $x$ because
$U$ acts freely on $P$. This shows that the stabilizer subgroup of $x$ in $N'\times U$ is finite. $\Box$

\lemma{
The principal $U$-bundle $P=Fr(\cN)$ over $Z$
has a $G\times U$-invariant connection one-form $\theta$
which vanishes on each tangent space to an $N'$-orbit in $P$. We shall use this one-form in our definition of the minimal-coupling symplectic form $\sigma^\#$.}
\thmlab\VanishingLemma
\proof
Fix a basis $\xi_i$ for $\gu$ and
a basis $\eta_j$ for $\gn'$. By the preceding lemma,
there exists one-forms $\theta^i,o^j$ on $P$ such that
$\iota_{\xi_i}\theta^j=\delta^j_i$, $\iota_{\eta_i}o^j=\delta^j_i$,
and $\iota_{\xi_i}o^j=0=\iota_{\eta_i}\theta^j$.
Then $\theta':=\sum\theta^i\otimes\xi^i$ is a connection one-form which vanishes on each
tangent space to an $N'$-orbit in $P$. Now average $\theta'$ over $G\times U$
and get a $G\times U$-invariant connection one-form $\theta$ with the same vanishing property.
$\Box$

Consider a hermitian vector bundle $E\ra M$,
and identify $E=Fr(E)\times_{U(n)}\C^n$. Then a tangent vector
of the shape $[*,0]\in TE=T(Fr(E))\times T\C^n)/U(n)$ is tangent to the
zero section $M\hra E$. We apply this to the following situation. Recall the diagram \NDiagram:
$$\matrix{
\cN^\# &\ra & \cN\cr
\pi_{\cN^\#}\da& &\da\pi_\cN\cr
\cV^*&\br\pi_{\cV^*}\over\ra& Z.
}$$
From the zero section $Z\hra\cV^*$, we get $\cN^\#|Z=\cN$ and 
the inclusion
$$
\cN=Fr(\cN)\times_U V=Fr(\cN^\#|Z)\times_U V\hra Fr(\cN^\#)\times_UV=\cN^\#.
$$
Since $Z\hra\cN$ as the zero section, we also have $Z\hra\cN^\#$ as a submanifold, and so
$T_zZ\subset T_z\cN\subset T_z\cN^\#$ for $z\in Z$.

\lemma{Consider a tangent vector of the shape $(y_1,0)\in T(Fr(\cN^\#)\times V)$
and its image $[y_1,0]\in T_z\cN^\#$ where $z\in Z$. If $[y_1,0]$ is tangent to $\cN$, then it
is tangent to $Z\hra\cN$.}
\thmlab\TangentToZ
\proof
By assumption $[y_1,0]\in T_z\cN$. Applying the observation above 
to the vector bundle $\cN\ra Z$ with typical fiber $V\equiv\C^n$,
we see that a tangent vector of the form $[y_1,0]\in T_z\cN$ is tangent to the zero section $Z$.
$\Box$

Extend the map $pr_F:Z\ra F$ to $\cN\supset Z$ by composing with 
the projection $\cN\ra Z$. Likewise extend it to $\cV^*$ and to $\cN^\#$.
Denote the extension of $pr_F$ by $pr_F^\#:\cN^\#\ra F\subset\gg^*$.

\lemma{Consider the principal bundle $P=Fr(\cN^\#)\ra\cV^*$, and  a tangent vector of the shape 
$[0,y_2]\in T(P\times_U V)=T\cN^\#$,
which is the image of tangent vector $(0,y_2)\in T(P\times V)$.
Then $\bra d(\xi\circ pr^\#_F),[0,y_2]\ket=0$ for any $\xi\in\gn$.}
\proof
By definition of $pr^\#_F$ as a composition of maps, we have
$$
d(\xi\circ pr^\#_F)=d(\xi\circ pr_F)\circ(\pi_{\cV^*})_*\circ(\pi_{\cN^\#})_*.
$$
Since $[0,y_2]$ is vertical it follows that $(\pi_{\cN^\#})_*[0,y_2]=0$. $\Box$

\lemma{Keep the same notations as in the preceding lemma. Consider $\xi\in\gn$ and a tangent vector of the shape
$[y_1,0]\in T(P\times_U V)=T\cN^\#$, 
which is the image of a tangent vector $(y_1,0)\in T(P\times V)$.
If $[y_1,0]$ is tangent to $\cN$, then
$\bra d(\xi\circ pr^\#_F),[y_1,0]\ket=\sigma^\#(X^\xi,[y_1,0])$.}
\proof
By Lemma \TangentToZ, a vector of the shape $[y_1,0]$ in $T\cN^\#$ being tangent to $\cN$
means that it is tangent to $Z\hra\cN$. 
By Lemma \RestrictionOfJ,  for $\tilde q=[dg,0,0]\times v\in Z$, we have 
$J_1(\tilde q)=pr_F(\tilde q)-v_2=v_1$, i.e.
$J_1=pr_F-v_2$ on $Z$. Since $\xi\circ v_2=0$ because $\xi\in\gn$ and $v_2\in\gd^*$,
it follows $d(\xi\circ pr^\#_F)=d(\xi\circ pr_F)=d(\xi\circ J_1)$
when evaluated on $[y_1,0]\in TZ$. Since $J_1$ is a $N'$-moment map
for the symplectic manifold ($\cN^\#,\sigma^\#)$, it follows that
$$
\bra d(\xi\circ J_1),[y_1,0]\ket=\sigma^\#(X^\xi,[y_1,0])
$$
which implies our assertion. $\Box$

\lemma{The restriction of $J_1$ to $\cN\hra\cN^\#$ coincides with
$pr_F^\#-v_2$. In fact we have $J_1:\cN=G\times_HV\times F\ra\gn^*$, $[g,\nu,v]\mapsto v_1$ where $v=(v_1,v_2)$.}
\thmlab\JOne
\proof
The second assertion follows from the first assertion.
By Lemma \RestrictionOfJ, $J_1=pr_F-v_2$
on each connected component of $Z$ ($G/H$ need not be connected).  
Since $\cN$ is a vector bundle over $Z$, the connected components
of $\cN$ corresponds 1-1 with the connected components in $Z$.
It suffices to show that the vector valued function $pr_F^\#-v_2$
satisfies $d(pr_F^\#-v_2)=dJ_1$ on $\cN$, i.e. that this equality holds as $\gn^*$-valued
function when contracted with each tangent vector $y\in T\cN$. Since $J_1$ is an $N'$-moment map
for $(\cN^\#,\sigma^\#)$, we have
$\bra d(\xi\circ J_1),y\ket=\bra \iota_\xi\sigma^\#,y\ket=\sigma^\#(X^\xi,y)$ for all $\xi\in\gn$.
It suffices to show that $pr_F^\#-v_2$ satisfies the same condition.
Since $\xi\circ v_2=0$ for $\xi\in\gn$ because $v_2\in\gd^*$, we only need to show that
$$
\bra d(\xi\circ pr_F^\#),y\ket=\sigma^\#(X^\xi,y),~~~\forall\xi\in\gn,~\forall y\in T\cN.
$$

As before, represent $y\in T\cN\subset T\cN^\#$ as the image $[y_1,y_2]$ of
tangent vector $(y_1,y_2)$ on $Fr(\cN^\#)\times V$ under the $U$-orbit map.
Since $y=[y_1,0]+[0,y_2]$, by the two preceding lemmas, it remains to show that
$$
\sigma^\#(X^\xi,[0,y_2])=0.
$$
Here is the calculations:
$$\eqalign{
&\sigma^\#(X^\xi,[0,y_2])\cr
&=(\sigma+\omega_V)([{\tilde X}^\xi,0,0],[0,\Psi_*(y_2),y_2]),\cr
&~~~~~~~\cN^\#=Fr(\cN^\#)\times_UV{\br i\over\cong}
(Fr(\cN^\#)\times\gu^*\times V)_0/U,~~[e,\nu]\mapsto[e,\Psi(\nu),\nu]\cr
&=\sigma([{\tilde X}^\xi,0],[0,\Psi_*(y_2)])+\omega_V(0,y_2),
~~~~\sigma~operates~on~1st~two~slots,~\omega_V~on~3rd~slot\cr
&=\mu((\tilde\pi_{\cN^\#})_*{\tilde X}^\xi,(\tilde\pi_{\cN^\#})_*0)-
d\bra pr_2,\tilde\theta\ket([{\tilde X}^\xi,0],[0,\Psi_*(y_2)]),~~~~
\sigma=(\tilde\pi_{\cN^\#})^*\mu-d\bra pr_2,\tilde\theta\ket\cr
&=-[{\tilde X}^\xi,0]'  \bra\bra pr_2,\tilde\theta\ket, [0,\Psi_*(y_2)]' \ket
+[0,\Psi_*(y_2)]'  \bra\bra pr_2,\tilde\theta\ket, [{\tilde X}^\xi,0]' \ket\cr
&~~~~+\bra\bra pr_2,\tilde\theta\ket, [ [{\tilde X}^\xi,0]', [0,\Psi_*(y_2)]' ]\ket.
}$$
For the last equality, we have use $d\alpha(X,Y)=X\alpha(Y)-Y\alpha(X)-\alpha([X,Y])$
for a one form $\alpha$ and vector fields $X,Y$. Here $[\cdots]'$ means
extending the tangent vector at a point to a tangent vector field in a neighborhood.
Consider the last three-term expression.
The first term is zero because $\tilde\theta$, being a one-form on $Fr(\cN^\#)$
operates on the first slot of $[0,\Psi_*(y_2)]'$. The third term is zero because
$[{\tilde X}^\xi,0]'$ can be viewed as a vector field on $Fr(\cN^\#)$, 
while $[0,\Psi_*(y_2)]'$ on $\gu^*$, hence they commute. Finally,
in the second term, the second factor is 
$\bra pr_2,\tilde\theta({\tilde X}^\xi)\ket$
But ${\tilde X}^\xi$ is a vector field generated by the $N'$-action on $Fr(\cN^\#)$
and $\tilde\theta$ is by construction vanishing on such vector fields, by Lemma \VanishingLemma.
It follows that this term is also zero. This completes the proof. $\Box$

\subsec{$T$-moment map near $X_F$}

In this subsection, we shall use the results we obtained about
the $G$-moment map near $Z\subset \cN$ to give a partial description of the $T$-moment map
after symplectically reducing $\cN^\#$ (near $Z$) with respect to $A$. As before, we denote by $J^\#$ the $G$-moment map on $\cN^\#$ corresponding to $J_G$ on $M$.

Recall that $J^\#=(J_1,J_2)$ takes values in $\gg^*=\gn^*\times\gd^*$.
When we reduce $\cN^\#$ with respect to $A$, we get
a semi-global form for $X_F$ in $X$ equipped with
a $T$-moment map. Since $A\subset D$, the $A$-reduction
does not affect $J_1$. On the $A$-reduced space,
$J_2$ induces a $(\gd/\ga)^*$ valued map $J_2'$ 
defined on $A$-orbits of the zero level set of $i^*\circ J_2$ where $i:A\hra D$.
The $T$-moment map on the $A$-reduced space is
then $J_F:=(J_1,J_2')$.
Note that $J_F$ is also $\gn^*\times\gb^*=\gt^*$-valued because
$\gd/\ga=\gh/\gk=\gb$ by \MainDiagram.
To describe $J_F$ explicitly, we will use the description of
$J_1|\cN$ in the preceding subsection.

The semi-global form gives a local description of a neighborhood
of $J^{-1}_G(F)$ in $M$ and the corresponding neighborhood of $J_T^{-1}(F)$ in $X$, in terms of
the data $H,K,D,V,\Phi$, corresponding to a chosen point $\tilde p\in J^{-1}_G(F)$. 
Let $\Psi$ be the $K=H\cap A$-moment map on $V$ with respect
to the linear action $K\subset H\ra U(V)$. 
 Consider the fiber \FiberOverXF.

\lemma{$A\backslash (Dg\times_H V)\cong tB\times_B V/K$.}
\proof
Since $D=A\cdot H$, we have $\pi(D)=\pi(A)\pi(H)=B$,
hence $\pi(Dg)=tB$. Thus we have a natural map
$$
\pi\times(-)/K:Dg\times V\thra tB\times V/K.
$$
This map is $H$-equivariant where $H$ acts on the domain by
$(dg,\nu)\mapsto(dgh^{-1},h\nu)$, and on the target via $\pi:H\ra B$.
Thus the map descends to a map $f:Dg\times_H V\thra tB\times_B V/K$.
It is straightforward to check that this map is 1-1. $\Box$

\lemma{Consider $\Psi:=i^*\circ\Phi$ the $K$-moment map on $V$ where $i:K\hra H$.
Then the $A$-zero level set in the fiber \FiberOverXF~ is
$$
(**)~~~~gD\times_H(0\times \Psi^{-1}(0))\times v.
$$
In particular, the $A$-reduced space is fiberwise $tB\times_B\Psi^{-1}(0)/K\times v$.}
\thmlab\FiberwiseReduction
\proof
In section 7 \SL, it was shown that the $A$-reduction on a symplectic bundle can be done fiberwise.
We apply this to the symplectic bundle $\cN^\#\ra X_F$
with fiber \FiberOverXF. This fiber has the $D$-moment map
$$
J_2([dg,\eta,\nu]\times v)=Ad^*(dg)(\eta+\Phi(\nu))+v_2.
$$
The $A$-zero level set on this fiber is defined by $i^*\circ J_2=0$.
Now $i^*:\gd^*\thra\ga^*$ is a $Ad^*(G)$-module map.
Moreover we have $i^*(v_2)=0$ because $v=(v_1,v_2)\in\gn^*\times\gb^*$
and $\gb^*=(\gd/\ga)^*$. So we get
$$
i^*\circ J_2([dg,\eta,\nu]\times v)=Ad^*(dg)(i^*\eta+\Psi(\nu)).
$$
This is zero iff $i^*\eta\in(\ga/\gk)^*$ and $\Psi(\nu)\in\gk^*$ are zero separately.
Note that $\eta\in (\gd/\gh)^*$ and $i:A/K\cong D/H$ is isomorphism (earlier lemma),
which means that $i^*\eta\in (\ga/\gk)^*$ is zero iff $\eta=0$.
This shows that on the fiber \FiberOverXF, the $A$-zero level set is exactly (**). 
Now performing reduction fiberwise, we get $(**)/A\cong tB\times_B\Psi^{-1}(0)/K\times v$.
$\Box$

\lemma{Let The $T$-moment map $J_F=(J_1,J_2')$ induced on
the $A$-reduced space $(i^*\circ J^\#)^{-1}(0)/A$ is given, on each fiber 
$tB\times_B\Psi^{-1}(0)/K\times v$ over $(tB,v)\in T/B\times F\subset X_F$, by
$$
J_F([tb,\nu]\times v)=v+\Theta(\nu)
$$
where $\Theta:\Psi^{-1}(0)/K\ra\gb^*=(\gh/\gk)^*$ 
is the $B$-moment map induced by $\Phi:V\ra\gh^*$
on the $K$-reduced space $\Psi^{-1}(0)/K$.}
\proof
By the preceding lemma, on each fiber the $A$-zero level set is (**),
having zero component along $(\gd/\gh)^*$. But this means that
this zero level set lies in $\cN\hra\cN^\#$. By 
Lemma \JOne, it follows that $J_1=pr_F^\#-v_2$ on the $A$-zero level set.
The same is true for the $A$-reduced space.

Now consider the second component of $J_F$. Again a point in the $A$-reduced space
has the shape $[[dg,0,\nu]]\times v$ (double bracket here means taking $A$-orbit), and we have
$$
J_2([[dg,0,\nu]]\times v)=J_2([tb,\nu]\times v)=\Phi(\nu)+v_2
$$
where $t=\pi(g)$ and $b=\pi(d)$.
Note that $Ad^*(dg)$ acts trivially on $\Phi(\nu)\in\gb^*$ because
$G$ acts trivially on $\gt$. It follows that
$$
J_F([tb,\nu]\times v)=v-v_2+\Phi(\nu)+v_2=v+\Phi(\nu).
$$
By definition $\Phi(\nu)$ is the value of the $B$-moment map on the $K$-orbit
of the vector $\nu\in V$. This completes the proof.  $\Box$



\corollary{Let $E\subset F$ be a compact set. There is closed ball $\B\subset\gb^*$ around zero, a $B$-invariant closed neighborhood $\cU\subset\Psi^{-1}(0)/K\subset V/K$ around zero, and a bijection $(W_{v_0}\times\B)\cap\Delta\ra \cU/B\times W_{v_0}$, $v+\beta\ra (B\cdot\nu,v)$ with
$$
\Theta(\nu)=\beta,~~~
J_F^{-1}(v+\beta)=T\times_B(B\cdot\nu)\times v.
$$
Moreover for any $\lambda\in[0,1]$, $v+\lambda^2\beta$ corresponds to $(B\cdot\lambda\nu,v)$.}
\thmlab\SemiGlobalII
\proof
The preimage of a small neighborhood of $v_0$ in $\Delta$ under $J_F$ is a small neighborhood of the $T$-orbit $J_F^{-1}(v_0)=T\times_B0\times v_0$ in $T\times_B\Psi^{-1}(0)/K\times F$, by the
preceding lemma. In particular for a small closed ball $\B\subset\gb^*$ around 0,
$J_F^{-1}(v+\B)$ is a closed $T$-invariant subset of $T\times_B\Psi^{-1}(0)/K\times v$, for each $v\in W_{v_0}$.
Thus there is a closed $B$-invariant neighborhood $\cU$ around 0 in $\Psi^{-1}(0)/K$ such that
$$
J_F^{-1}(v+\B)=T\times_B\cU\times v.
$$
By the Orbit Theorem, for each $v+\beta\in(v+\B)\cap\Delta$, there is a unique $B$-orbit $B\cdot\nu$
in $\cU$ such that $J_F^{-1}(v+\beta)=T\times_B(B\cdot\nu)\times v$.
Then $v+\beta=v+\Theta(\nu)$ follows from the preceding lemma. This proves the first assertion.
The second assertion follows from that $\Theta(\lambda\nu)=\lambda^2\Theta(\nu)$.
$\Box$

In this subsection, we will prove Theorem \TheoremIII.
Thus we are given two exact sequences $1\ra A(k)\ra G(k)\ra T\ra 1$, $k=1,2$, of Lie groups
and symplectic manifolds $(M(k),\omega(k))$ such that the respective $A(k)$-reductions $X(k)$
have identical moment polytope $\Delta=J_T(1)(X(1))=J_T(2)(X(2))$ in $\gt^*$. We would like to construct a $T$-equivariant homeomorphism $\varphi:X(1)\ra X(2)$ such that $J_T(1)=J_T(2)\circ\varphi$, under the assumption of Theorem \TheoremIII.

\bs
{\it Step 0.} We first construct a covering of $\Delta$ by compact sets as follows.
Fix a vertex $E$ of $\Delta$. For $\epsilon>0$, let $\B(\epsilon)\subset\gt^*$ be
the closed ball centered at 0 of radius $\epsilon$. We identify $E\times\B(\epsilon)$
with the $\epsilon$-ball centered at $E$. By Corollary 
\SemiGlobalII, we can choose $\epsilon$ so that a semi-global form
corresponding to $E$ is valid in $J_T^{-1}(E\times\B(\epsilon))$. Likewise do the same
for each vertex, and shrink the $\epsilon$ if necessary so that the
$E\times\B(\epsilon)$ do not overlap.  
We also arrange that these properties continue holds when the $\epsilon$ are perturbed slightly. 
Next fix a 1-face and consider its interior $F$.
Cut off both ends of $F$ slightly to get a compact set $E\subset F$.
For $\epsilon>0$, let $\B(\epsilon)\subset\gb^*$ ($\gb^*$ depends on $F\subset v+\gn^*$)
be the closed ball centered at $0$ of radius $\epsilon$. Again we identify
$E\times\B(\epsilon)$ with a compact tubular neighborhood of $E$ obtained
by thickening $E$ in the normal directions (i.e. along $\gb^*$) by $\epsilon$.
The identification is $E\times\B(\epsilon)\ni(v,\beta)\leftrightarrow v+\beta\in\gt^*$.
Choose $\epsilon$ so that a semi-global form
corresponding to $F$ is valid in $J_T^{-1}(E\times\B(\epsilon))$. Likewise do the same
for each 1-face, and shrink the $\epsilon$ if necessary so that the
$E\times\B(\epsilon)$ corresponding to the 1-faces do not overlap. 
Likewise do the same for all 2-faces, 3-faces, ..., $(n-1)$-faces.
The result is a covering of $\Delta$ by compact neighborhoods of the shape $E\times\B(\epsilon)$,
one for each proper face of $\Delta$, with the properties that those neighborhoods supported
on faces of the same dimension do not overlap, and that a semi-global form
is valid in each $J_T^{-1}(E\times\B(\epsilon))$. Again, we arrange the covering so that
all these properties continue to hold when the $\epsilon$ are perturbed slightly.

Next, we shall use the semi-global form for each face $F$ to define a global homeomorphism.

\lemma{Let $S$ be the set of $K$-orbits $\nu\in\Psi^{-1}(0)/K$
such that $\Theta(\nu)\in\partial\B(\epsilon)$, and put
$S_\leq:=\{\lambda\nu|\lambda\in[0,1],~\nu\in S\}$. Then
\item{i.} $S$ is a subset of quadric hypersurface in the stratified space $V/K$.
\item{ii.}  $J_F^{-1}(E\times\partial\B(\epsilon))=T\times_B S\times E$.
\item{iii.} Every point in $J_F^{-1}(E\times\B(\epsilon))$
away from the zero section $X_F$ has the shape $[t,\lambda\nu]\times v$
for a unique $\lambda\in(0,1]$ and $[t,\nu]\in T\times_B S$. 
\item{iv.} $J^{-1}_F(E\times\B(\epsilon))=T\times_BS_\leq\times E$.
\item{v.} $(E\times\B(\epsilon))\cap\Delta\cong S_\leq/B\times E,~(v,\beta)\mapsto (B\cdot\nu,v)$.
\item{vi.} There is a dense subset $S^\circ\subset S$ on which $B$ acts
freely such that 
$$
J_F^{-1}((E\times\partial\B(\epsilon))\cap\Delta^\circ)=T\times_B S^\circ\times E.
$$
}\thmlab\Filling
\proof
Part i. follows from that $\Psi(\nu)=0$ is a quadratic equation on $V$.
Part ii.-iv. follow from Corollary \SemiGlobalII. 
Part v. follows from taking $T$-orbit spaces on both
sides of iv. and applying the Orbit Lemma on the left hand side.

By a corollary to the Face Lemma, $T$ acts freely on
$J_F^{-1}((E\times\partial\B(\epsilon))\cap\Delta^\circ)$.
By part ii., we have $J_F^{-1}((E\times\partial\B(\epsilon))\cap\Delta^\circ)=T\times_B S^\circ\times E$
for some $S^\circ\subset S$. Since the $T$-action on the left hand side is free,
the $B$-action on $S^\circ$ is also free. In fact $S^\circ$ must be the
full subset of $S$ on which $B$ acts freely.
Since $J_F$ is an open mapping it follows that
the closure of $J_F^{-1}((E\times\partial\B(\epsilon))\cap\Delta^\circ)$
coincides with $J_F^{-1}(E\times\partial\B(\epsilon))=T\times_B S\times E$.
This means that $S$ is the closure of $S^\circ$. $\Box$

\bs
{\it Step 1.} We now begin dealing with two reduced Delzant $T$-spaces $X(k)$ with the same moment polytope $\Delta$, as before. Thus each has its own moment map $J_T(k)$, and for each face $F$ of $\Delta$, $X(k)$ has its own moment map $J_F(k)$ defined on a semi-global form corresponding to a compact neighborhood $E\subset F$ constructed in Step 0. Corresponding to this is a stabilizer subgroup $B(k)$ in $T$, and quadric $S(k)$ as in Lemma \Filling.
We begin with a fixed $(n-1)$-face $F$. Then we have the following $T$-equivariant
$J_T$-compatible commutative diagram
$$\matrix{
T\times\Delta^\circ&=&T\times\Delta^\circ\cr
\sim\ua& &\ua\sim\cr
J_T(1)^{-1}(\Delta^\circ)&\ra&J_T(2)^{-1}(\Delta^\circ)\cr
\cup & & \cup\cr
J_T(1)^{-1}(E\times\partial\B(\epsilon))&\ra &J_T(2)^{-1}(E\times\partial\B(\epsilon))\cr
\cap & & \cap\cr
J_T(1)^{-1}(E\times\B(\epsilon)) &\cdots\ra &J_T(2)^{-1}(E\times\B(\epsilon)).
}$$
The first horizontal arrow is an isomorphism induced by the equality.
The second horizontal arrow is the restriction of the first one
to one side of the boundary of the neighborhoods.
The third horizontal arrow is an extension, which we want to construct, of the second arrow
to the interior of the neighborhoods. 
We shall do so by using Lemma \Filling. 
More precisely we want to complete the diagram $T$-equivariantly, $J_T$-compatibly, 
obtained from the bottom rectangle with $J_T(k)$ replaced by $J_F(k)$.

By Lemma \Filling, the problem becomes completing the diagram
\eqn\Extension{\matrix{
T\times_{B(1)}S(1)\times E &\br\varphi_F\over \ra& T\times_{B(2)} S(2)\times E\cr
\cap & & \cap\cr
T\times_{B(1)}S(1)_\leq\times E & \cdots\ra& T\times_{B(2)} S(2)_\leq\times E.
}}
By assumption of Theorem \TheoremIII, the $B(k)$ are connected.
Since the $B(k)$ have the same Lie algebra $\gb\subset\gt$, it follows that $B(1)=B(2)$
which we will denote by $B$.

\lemma{$\varphi_F$ preserves the factor $E\subset F$.}
\proof
Since $\varphi_F$ is the restriction of a $J$-compatible map,
if we write $\varphi_F([t_1,\nu_1]\times v_1)=[t_2,\nu_2]\times v_2$ and apply $J_F(2)$,
we get $v_1+\Theta(1)(\nu_1)=v_2+\Theta(2)(\nu_2)\in E\times\B(\epsilon)$.
Projecting into the $E$ factor, we get $v_1=v_2$. $\Box$

Since $dim~F=n-1$, the ball $\B(\epsilon)$ is an interval and $(v\times\partial\B(\epsilon))\cap\Delta$ is a single point in the interior $\Delta^\circ$, for each $v\in E$. It follows from Lemma \Filling~that
each $S(k)$ is a single free $B$-orbit, say $B\cdot\nu_k$. 
Hence there exists a unique continuous function $f:E\ra T$ such that
$$
\varphi_F([e,\nu_1]\times v)=[f(v),\nu_2]\times v.
$$
Extend this to $T\times_B S(1)_\leq\times E\ra T\times_B S(2)_\leq\times E$ by
$$
\tilde\varphi_F:[t,\lambda b\cdot\nu_1]\times v\mapsto [tf(v),\lambda b\cdot\nu_2]\times v.
$$
By Lemma \Filling, this is well-defined for all $t\in T, b\in B$, and $\lambda\in[0,1]$; it is also
manifestly $T$-equivariant. It is clear that this is a homemorphism.
Using the quadratic scaling property of $\Theta(k)$, it is easy
to check that the extension is compatible with the moment maps $J_F(k)$.
Hence $\tilde \varphi_F$ gives a well-defined map 
$$
J_T(1)^{-1}(E\times\B(\epsilon))\ra J_T(2)^{-1}(E\times\B(\epsilon)).
$$

Thus $\varphi_{\Delta^\circ}$ together with the $\tilde\varphi_F$ constructed above for $dim~F=n-1$, defines a $T$-equivariant homeomorphism 
$$
\varphi: J_T(1)^{-1}(\Delta_{n-1})\ra J_T(2)^{-1}(\Delta_{n-1})
$$ 
which is compatible with the moment maps $J_T(k)$. Here $\Delta_p\subset\Delta$ is the closure of
the complement of the union of the compact neighborhoods $E\subset F$ for $dim~F<p$.

{\it Step 2.} $dim~F=n-2$ or lower. We proceed by induction. Thus suppose that we have defined a
$T$-equivariant homeomorphism 
$$
\varphi: J_T(1)^{-1}(\Delta_p)\ra J_T(2)^{-1}(\Delta_p)
$$ 
which is compatible with the moment maps $J_T(k)$ for a given $p\leq n-1$. Then for each $F$ of dimension $dim~F=p-1$, the top row $\varphi_F$ in the diagram \Extension~ is a well-defined $T$-equivariant homeomorphism compatible with the moment maps. We want to complete that diagram, i.e. to extend the $\varphi_F$ to a $T$-equivariant homemorphism $\tilde\varphi_F$ which is compatible with the moment maps. 

Since $B\subset T$ is a assumed connected, we have a splitting
$$
T=N\times B.
$$
It follows that
$$
T\times_B S(k)_\leq\times E=N\times S(k)_\leq\times E
$$
as $T$-spaces, where $T$ acts on the first two factors of the right hand side via $T=N\times B$.
Thus for the top row of \Extension,
there exist unique continuous functions $f,g$ such that
$$
\varphi_F: N\times S(1)\times E\ra N\times S(2)\times E,
~~(n,\nu)\times v\mapsto(n\cdot f(\nu,v),g(\nu,v))\times v.
$$
Since $\varphi_F$ is also $B$-equivariant, it follows that $f:S(1)\times E\ra N$ is $B$-invariant,
and $g:S(1)\times E\ra S(2)$ is $B$-equivariant.

\lemma{There exists a homotopy
$\phi:[0,1]\times S(1)\times E\ra N$ such that  $\phi(1,\nu,v)=f(\nu,v)$, $\phi(0,\nu,v)=e$, and that the map $\tilde \varphi_F:N\times S(1)_\leq\times E\ra N\times S(2)_\leq\times E$ with
$$
\tilde\varphi_F(n,\lambda\nu,v)= (n\cdot\phi(\lambda,\nu,v),\lambda g(\nu,v))\times v
$$
is a homeomorphism which is the identity map when $\lambda=0$. Moreover $\tilde\varphi_F$ is
$T$-equivariant and compatible with the moment maps $J_F(k)$.}
\proof
By Lemma \Filling$v$, $S(1)/B$ is homeomorphic to a closed ball. In particular it is contractible, i.e. there is a continuous map
$$
r:[0,1]\times S(1)/B\ra S(1)/B
$$
such that $r(1,-)=id$ and $r(0,-)=pt.$ Put $\phi(\lambda,\nu,v)=f(\nu',v)$ where $\nu'\in r(\lambda,B\cdot\nu)$. Since $f$ is $B$-invariant, this is independent of the choice of $\nu'$. 
If the point $r(1,B\cdot\nu)=pt$ does not get mapped to $e\in N$ under $f$,
then multiply $\phi$ by a curve $c:[0,1]\ra N$ connecting $e$ to $f(pt)$.
Thus $\phi$ satisfies the boundary conditions we seek.
That $\tilde\varphi_F$ is a homeomorphism follows easily from
Lemma \Filling. 

Finally, $T=N\times B$-equivariance of $\tilde\varphi_F$ follows readily from the $B$-equivariance of $g$, the $B$-invariance of $f$. That $\tilde\varphi_F$ is compatible with the $J_F(k)$ follows from
the preceding lemma and that the maps $\Theta(k)$ are quadratic. $\Box$

Since the compact neighborhoods $E$ we attach to the faces $F$ of a given dimension $p-1$
do not overlap by construction, the same procedure can be performed separately to construct the desired extension $\tilde\varphi_F$ of $\varphi_F$, for each such $F$. Together with our inductive hypothesis,
we have now extended $\varphi$ above to a $T$-equivariant homeomorphism 
$$
\varphi: J_T(1)^{-1}(\Delta_{p-1})\ra J_T(2)^{-1}(\Delta_{p-1})
$$ 
which is compatible with the moment maps $J_T(k)$. This completes the proof of Theorem \TheoremIII. $\Box$ 

\footatend\vfill\supereject\immediate\closeout\rfile\writestoppt
\baselineskip=14pt\centerline{{\bf References}}\bigskip{\frenchspacing%
\parindent=20pt\escapechar=` \input refs.tmp\vfill
\eject}\nonfrenchspacing

\bs
\item{} Bong H. Lian, Department of Mathematics, National University of Singapore, 2 Science Drive 2, Singapore 117543. On leave of absence from Department of Mathematics, Brandeis University, Waltham MA 02454.
\bs
\item{} Bailin Song, Department of Mathematics, Brandeis University, Waltham MA 02454.

\end

\subsec{Construction of $T$-homeomorphism}

In this subsection, we will prove Theorem \TheoremIII.
Thus we are given two exact sequences $1\ra A(k)\ra G(k)\ra T\ra 1$, $k=1,2$, of Lie groups
and symplectic manifolds $(M(k),\omega(k))$ such that the respective $A(k)$-reductions $X(k)$
have identical moment polytope $\Delta=J_T(1)(X(1))=J_T(2)(X(2))$ in $\gt^*$. We would like to construct a $T$-equivariant homeomorphism $\varphi:X(1)\ra X(2)$ such that $J_T(1)=J_T(2)\circ\varphi$, under the assumption of Theorem \TheoremIII.

\bs
{\it Step 0.} We first construct a covering of $\Delta$ by compact sets as follows.
Fix a vertex $E$ of $\Delta$. For $\epsilon>0$, let $\B(\epsilon)\subset\gt^*$ be
the closed ball centered at 0 of radius $\epsilon$. We identify $E\times\B(\epsilon)$
with the $\epsilon$-ball centered at $E$. By Corollary 
\SemiGlobalII, we can choose $\epsilon$ so that a semi-global form
corresponding to $E$ is valid in $J_T^{-1}(E\times\B(\epsilon))$. Likewise do the same
for each vertex, and shrink the $\epsilon$ if necessary so that the
$E\times\B(\epsilon)$ do not overlap.  
We also arrange that these properties continue holds when the $\epsilon$ are perturbed slightly. 
Next fix a 1-face and consider its interior $F$.
Cut off both ends of $F$ slightly to get a compact set $E\subset F$.
For $\epsilon>0$, let $\B(\epsilon)\subset\gb^*$ ($\gb^*$ depends on $F\subset v+\gn^*$)
be the closed ball centered at $0$ of radius $\epsilon$. Again we identify
$E\times\B(\epsilon)$ with a compact tubular neighborhood of $E$ obtained
by thickening $E$ in the normal directions (i.e. along $\gb^*$) by $\epsilon$.
The identification is $E\times\B(\epsilon)\ni(v,\beta)\leftrightarrow v+\beta\in\gt^*$.
Choose $\epsilon$ so that a semi-global form
corresponding to $F$ is valid in $J_T^{-1}(E\times\B(\epsilon))$. Likewise do the same
for each 1-face, and shrink the $\epsilon$ if necessary so that the
$E\times\B(\epsilon)$ corresponding to the 1-faces do not overlap. 
Likewise do the same for all 2-faces, 3-faces, ..., $(n-1)$-faces.
The result is a covering of $\Delta$ by compact neighborhoods of the shape $E\times\B(\epsilon)$,
one for each proper face of $\Delta$, with the properties that those neighborhoods supported
on faces of the same dimension do not overlap, and that a semi-global form
is valid in each $J_T^{-1}(E\times\B(\epsilon))$. Again, we arrange the covering so that
all these properties continue to hold when the $\epsilon$ are perturbed slightly.

Next, we shall use the local normal form to define our homeomorphism.

\lemma{Let $S$ be the set of $K$-orbits $\nu\in\Psi^{-1}(0)/K$
such that $\Theta(\nu)\in\partial\B(\epsilon)$, and put
$S_\leq:=\{\lambda\nu|\lambda\in[0,1],~\nu\in S\}$. Then
\item{i.} $S$ is a subset of quadric hypersurface in the stratified space $V/K$.
\item{ii.}  $J_F^{-1}(W_{v_0}\times\partial\B(\epsilon))=T\times_B S\times W_{v_0}$.
\item{iii.} Every point in $J_F^{-1}(W_{v_0}\times\B(\epsilon))$
away from the zero section $X_F$ has the shape $[t,\lambda\nu]\times v$
for a unique $\lambda\in(0,1]$ and $[t,\nu]\in T\times_B S$. 
\item{iv.} $J^{-1}_F(W_{v_0}\times\B(\epsilon))=T\times_BS_\leq\times W_{v_0}$.
\item{v.} $(W_{v_0}\times\B(\epsilon))\cap\Delta\cong S_\leq/B\times W_{v_0},~(v,\beta)\mapsto (B\cdot\nu,v)$.
\item{vi.} There is a dense subset $S^\circ\subset S$ on which $B$ acts
freely such that 
$$
J_F^{-1}((W_{v_0}\times\partial\B(\epsilon))\cap\Delta^\circ)=T\times_B S^\circ\times W_{v_0}.
$$
}\thmlab\Filling
\proof
Part i. follows from that $\Psi(\nu)=0$ is a quadratic equation on $V$.
Part ii.-iv. follow from Corollary \SemiGlobalII. 
Part v. follows from taking $T$-orbit spaces on both
sides of iv. and applying the Orbit Lemma on the left hand side.

By a corollary to the Face Lemma, $T$ acts freely on
$J_F^{-1}((W_{v_0}\times\partial\B(\epsilon))\cap\Delta^\circ)$.
By part ii., we have $J_F^{-1}((W_{v_0}\times\partial\B(\epsilon))\cap\Delta^\circ)=T\times_B S^\circ\times E$
for some $S^\circ\subset S$. Since the $T$-action on the left hand side is free,
the $B$-action on $S^\circ$ is also free. In fact $S^\circ$ must be the
full subset of $S$ on which $B$ acts freely.
Since $J_F$ is an open mapping it follows that
the closure of $J_F^{-1}((W_{v_0}\times\partial\B(\epsilon))\cap\Delta^\circ)$
coincides with $J_F^{-1}(W_{v_0}\times\partial\B(\epsilon))=T\times_B S\times W_{v_0}$.
This means that $S$ is the closure of $S^\circ$. $\Box$

\bs
{\it Step 1.} We now begin dealing with two reduced Delzant $T$-spaces $X(k)$ with the same moment polytope $\Delta$, as before. Thus each has its own moment map $J_T(k)$ and local normal form $J_F(k)$ for each $F$ and some neighborhood $W_{v_0}\subset F$, its own stabilizer subgroup $B(k)$ in $T$, and quadric $S(k)$ as in Lemma \Filling.
Fix an $(n-1)$-face $F$. Then we have the following $T$-equivariant
$J_T$-compatible commutative diagram
$$\matrix{
T\times\Delta^\circ&=&T\times\Delta^\circ\cr
\sim\ua& &\ua\sim\cr
J_T(1)^{-1}(\Delta^\circ)&\ra&J_T(2)^{-1}(\Delta^\circ)\cr
\cup & & \cup\cr
J_T(1)^{-1}(E\times\partial\B(\epsilon))&\ra &J_T(2)^{-1}(E\times\partial\B(\epsilon))\cr
\cap & & \cap\cr
J_T(1)^{-1}(E\times\B(\epsilon)) &\cdots\ra &J_T(2)^{-1}(E\times\B(\epsilon)).
}$$
The first horizontal arrow is an isomorphism induced by the equality.
The second horizontal arrow is the restriction of the first one
to one side of the boundary of the neighborhoods.
The third horizontal arrow is an extension, which we want to construct, of the second arrow
to the interior of the neighborhoods. 
We shall do so by first define the extension locally, i.e. over each neighborhood $W_{v_0}$, by using the local normal form valid there. Then it will be clear that the resulting extension is well-defined, i.e. independent of the choice of the local normal form. 
More precisely we want to complete the diagram $T$-equivariantly, $J_T$-compatibly, 
obtained from the bottom rectangle with $J_T(k)$ replaced by $J_F(k)$.

By Lemma \Filling, the problem becomes completing the diagram
\eqn\Extension{\matrix{
T\times_{B(1)}S(1)\times W_{v_0} &\br\varphi_F\over \ra& T\times_{B(2)} S(2)\times W_{v_0}\cr
\cap & & \cap\cr
T\times_{B(1)}S(1)_\leq\times W_{v_0} & \cdots\ra& T\times_{B(2)} S(2)_\leq\times W_{v_0}.
}}
Note that the $S(k)$ depends on the choice of the local normal form over $W_{v_0}$ for each of the spaces $X(k)$. It will be clear the extension \Extension~ is compatible with a change of local normal form. By assumption of Theorem \TheoremIII, the $B(k)$ are connected.
Since the $B(k)$ have the same Lie algebra $\gb\subset\gt$, it follows that $B(1)=B(2)$
which we will denote by $B$.

\lemma{$\varphi_F$ preserves the factor $W_{v_0}\subset F$.}
\proof
Since $\varphi_F$ is the restriction of a $J$-compatible map,
if we write $\varphi_F([t_1,\nu_1]\times v_1)=[t_2,\nu_2]\times v_2$ and apply $J_F(2)$,
we get $v_1+\Theta(1)(\nu_1)=v_2+\Theta(2)(\nu_2)\in W_{v_0}\times\B(\epsilon)$.
Projecting into the $W_{v_0}$ factor, we get $v_1=v_2$. $\Box$

Since $dim~F=n-1$, the ball $\B(\epsilon)$ is an interval and $(v\times\partial\B(\epsilon))\cap\Delta$ is a single point in the interior $\Delta^\circ$, for each $v\in E$. It follows from Lemma \Filling~that
each $S(k)$ is a single free $B$-orbit, say $B\cdot\nu_k$. 
Hence there exists a unique smooth function $f:W_{v_0}\ra T$ (depending on the local normal form) such that
$$
\varphi_F([e,\nu_1]\times v)=[f(v),\nu_2]\times v.
$$
Extend this to $T\times_B S(1)_\leq\times W_{v_0}\ra T\times_B S(2)_\leq\times W_{v_0}$ by
$$
\tilde\varphi_F:[t,\lambda b\cdot\nu_1]\times v\mapsto [tf(v),\lambda b\cdot\nu_2]\times v.
$$
By Lemma \Filling, this is well-defined for all $t\in T, b\in B$, and $\lambda\in[0,1]$; it is also
manifestly $T$-equivariant. It is clear that this is a homemorphism.
Using the quadratic scaling property of $\Theta(k)$, it is easy
to check that the extension is compatible with the moment maps $J_F(k)$.
Finally, note that any other choice of local normal form of a $T$-orbit is, by construction, $T$-equivariantly symplectomorphic to the any other choice. The $S(k), S(k)_\leq$ for different local normal forms are also $B$-equivariantly equivalent under the same symplectomorphism. It is clear that the map $\tilde \varphi_F$ defined above is compatible with such a symplectomorphism. Hence $\tilde \varphi_F$ gives a well-defined map 
$$
J_T(1)^{-1}(E\times\B(\epsilon))\ra J_T(2)^{-1}(E\times\B(\epsilon)).
$$

Thus $\varphi_{\Delta^\circ}$ together with the $\tilde\varphi_F$ constructed above for $dim~F=n-1$, defines a $T$-equivariant homeomorphism 
$$
\varphi: J_T(1)^{-1}(\Delta_{n-1})\ra J_T(2)^{-1}(\Delta_{n-1})
$$ 
which is compatible with the moment maps $J_T(k)$. Here $\Delta_p\subset\Delta$ is the closure of
the complement of the union of the compact neighborhoods $E\subset F$ for $dim~F<p$.

\bs
{\it Step 2.} $dim~F=n-2$ or lower. We proceed by induction. Thus suppose that we have defined a
$T$-equivariant homeomorphism 
$$
\varphi: J_T(1)^{-1}(\Delta_p)\ra J_T(2)^{-1}(\Delta_p)
$$ 
which is compatible with the moment maps $J_T(k)$ for a given $p\leq n-1$. Then for each $F$ of dimension $dim~F=p-1$, the top row $\varphi_F$ in the diagram \Extension~ is a well-defined $T$-equivariant homeomorphism compatible with the moment maps. We want to complete that diagram, i.e. to extend the $\varphi_F$ to a $T$-equivariant homemorphism $\tilde\varphi_F$ which is compatible with the moment maps, and see that the extension is independent of the choice of local normal form.

Since $B\subset T$ is a assumed connected, we have a splitting
$$
T=N\times B.
$$
It follows that
$$
T\times_B S(k)_\leq\times W_{v_0}=N\times S(k)_\leq\times W_{v_0}
$$
as $T$-spaces, where $T$ acts on the first two factors of the right hand side via $T=N\times B$.
Thus for the top row of \Extension,
there exist {\it unique} smooth functions $f,g$ such that
$$
\varphi_F: N\times S(1)\times W_{v_0}\ra N\times S(2)\times W_{v_0},
~~(n,\nu)\times v\mapsto(n\cdot f(\nu,v),g(\nu,v))\times v.
$$
Since $\varphi_F$ is also $B$-equivariant, it follows that $f:S(1)\times W_{v_0}\ra N$ is $B$-invariant,
and $g:S(1)\times W_{v_0}\ra S(2)$ is $B$-equivariant.

\lemma{There exists a smooth homotopy
$\phi:[0,1]\times S(1)\times W_{v_0}\ra N$ such that  $\phi(1,\nu,v)=f(\nu,v)$, $\phi(0,\nu,v)=e$, and that the map $\tilde \varphi_F:N\times S(1)_\leq\times W_{v_0}\ra N\times S(2)_\leq\times W_{v_0}$ with
$$
\tilde\varphi_F(n,\lambda\nu,v)= (n\cdot\phi(\lambda,\nu,v),\lambda g(\nu,v))\times v
$$
is a homeomorphism which is the identity map when $\lambda=0$. Moreover $\tilde\varphi_F$ is
$T$-equivariant and compatible with the moment maps $J_F(k)$.}
\proof
By Lemma \Filling$v$, $S(1)/B$ is canonically diffeomorphic to the closed ball $(v_0\times\partial B(\epsilon))\cap\Delta$. In particular it is smoothly contractible, i.e. there is a smooth map
$$
r:[0,1]\times S(1)/B\ra S(1)/B
$$
such that $r(1,-)=id$ and $r(0,-)=pt.$ Put $\phi(\lambda,\nu,v)=f(\nu',v)$ where $\nu'\in r(\lambda,B\cdot\nu)$. Since $f$ is $B$-invariant, this is independent of the choice of $\nu'$. Since $S(1)/B$ is canonically diffeomorphic to $(v_0\times\partial B(\epsilon))\cap\Delta$, $r$ can be made compatible with any change of local normal form.
If the point $r(1,B\cdot\nu)=pt$ does not get mapped to $e\in N$ under $f$,
then multiply $\phi$ by a curve $c:[0,1]\ra N$ connecting $e$ to $f(pt)$.
Thus $\phi$ satisfies the boundary conditions we seek.
That $\tilde\varphi_F$ is a homeomorphism follows easily from
Lemma \Filling. As in the codimension 1 case, the map is also compatible with any change of local normal form.

Finally, $T=N\times B$-equivariance of $\tilde\varphi_F$ follows readily from the $B$-equivariance of $g$, the $B$-invariance of $f$. That $\tilde\varphi_F$ is compatible with the $J_F(k)$ follows from
the preceding lemma and that the maps $\Theta(k)$ are quadratic. $\Box$

Since the compact neighborhoods $E$ we attach to the faces $F$ of a given dimension $p-1$
do not overlap by construction, the same procedure can be performed separately to construct the desired extension $\tilde\varphi_F$ of $\varphi_F$, for each such $F$. Together with our inductive hypothesis,
we have now extended $\varphi$ above to a $T$-equivariant homeomorphism 
$$
\varphi: J_T(1)^{-1}(\Delta_{p-1})\ra J_T(2)^{-1}(\Delta_{p-1})
$$ 
which is compatible with the moment maps $J_T(k)$. This completes the proof of Theorem \TheoremIII. $\Box$

\np
\footatend\vfill\supereject\immediate\closeout\rfile\writestoppt
\baselineskip=14pt\centerline{{\bf References}}\bigskip{\frenchspacing%
\parindent=20pt\escapechar=` \input refs.tmp\vfill
\eject}\nonfrenchspacing

\bs
\item{} Bong H. Lian, Department of Mathematics, National University of Singapore, 2 Science Drive 2, Singapore 117543. On leave of absence from Department of Mathematics, Brandeis University, Waltham MA 02454.
\bs
\item{} Bailin Song, Department of Mathematics, Brandeis University, Waltham MA 02454.

\end

%% file: abbrev.tex

%
%
%
\def\unredoffs{} \def\redoffs{\voffset=-.31truein\hoffset=-.59truein}
\def\speclscape{\special{ps: landscape}}
%
%
%
%
\newbox\leftpage \newdimen\fullhsize \newdimen\hstitle \newdimen\hsbody
\tolerance=1000\hfuzz=2pt
\catcode`\@=11 
%
\ifx\answ\bigans\message{(This will come out unreduced.}
\magnification=1200\unredoffs\baselineskip=16pt plus 2pt minus 1pt
\hsbody=\hsize \hstitle=\hsize 
\else\message{(This will be reduced.} \let\l@r=L
\magnification=1000\baselineskip=16pt plus 2pt minus 1pt \vsize=7truein
\redoffs \hstitle=8truein\hsbody=4.75truein\fullhsize=10truein\hsize=\hsbody
\output={\ifnum\pageno=0 
  \shipout\vbox{\speclscape{\hsize\fullhsize\makeheadline}
    \hbox to \fullhsize{\hfill\pagebody\hfill}}\advancepageno
  \else
  \almostshipout{\leftline{\vbox{\pagebody\makefootline}}}\advancepageno
  \fi}
\def\almostshipout#1{\if L\l@r \count1=1 \message{[\the\count0.\the\count1]}
      \global\setbox\leftpage=#1 \global\let\l@r=R
 \else \count1=2
  \shipout\vbox{\speclscape{\hsize\fullhsize\makeheadline}
      \hbox to\fullhsize{\box\leftpage\hfil#1}}  \global\let\l@r=L\fi}
\fi
%
\newcount\yearltd\yearltd=\year\advance\yearltd by -1900

%

%
%

\def\draftmode{\message{ DRAFTMODE }\def\draftdate{{\rm preliminary draft:
\number\month/\number\day/\number\yearltd\ \ \hourmin}}%
\headline={\hfil\draftdate}\writelabels\baselineskip=20pt plus 2pt minus 2pt
 {\count255=\time\divide\count255 by 60 \xdef\hourmin{\number\count255}
  \multiply\count255 by-60\advance\count255 by\time
  \xdef\hourmin{\hourmin:\ifnum\count255<10 0\fi\the\count255}}}
\def\nolabels{\def\wrlabeL##1{}\def\eqlabeL##1{}\def\reflabeL##1{}}
\def\writelabels{\def\wrlabeL##1{\leavevmode\vadjust{\rlap{\smash%
{\line{{\escapechar=` \hfill\rlap{\sevenrm\hskip.03in\string##1}}}}}}}%
\def\eqlabeL##1{{\escapechar-1\rlap{\sevenrm\hskip.05in\string##1}}}%
\def\reflabeL##1{\noexpand\llap{\noexpand\sevenrm\string\string\string##1}}}
\nolabels
%
\global\newcount\secno \global\secno=0
\global\newcount\meqno \global\meqno=1
\def\newsec#1{\global\advance\secno by1\message{(\the\secno. #1)}
\global\subsecno=0\eqnres@t\noindent{\bf\the\secno. #1}
\writetoca{{\secsym} {#1}}\par\nobreak\medskip\nobreak}
\def\eqnres@t{\xdef\secsym{\the\secno.}\global\meqno=1\bigbreak\bigskip}
\def\sequentialequations{\def\eqnres@t{\bigbreak}}\xdef\secsym{}
\global\newcount\subsecno \global\subsecno=0
\def\subsec#1{\global\advance\subsecno by1\message{(\secsym\the\subsecno. #1)}
\ifnum\lastpenalty>9000\else\bigbreak\fi
\noindent{\it\secsym\the\subsecno. #1}\writetoca{\string\quad
{\secsym\the\subsecno.} {#1}}\par\nobreak\medskip\nobreak}
\def\appendix#1#2{\global\meqno=1\global\subsecno=0\xdef\secsym{\hbox{#1.}}
\bigbreak\bigskip\noindent{\bf Appendix #1. #2}\message{(#1. #2)}
\writetoca{Appendix {#1.} {#2}}\par\nobreak\medskip\nobreak}
%
%
\def\eqnn#1{\xdef #1{(\secsym\the\meqno)}\writedef{#1\leftbracket#1}%
\global\advance\meqno by1\wrlabeL#1}
\def\eqna#1{\xdef #1##1{\hbox{$(\secsym\the\meqno##1)$}}
\writedef{#1\numbersign1\leftbracket#1{\numbersign1}}%
\global\advance\meqno by1\wrlabeL{#1$\{\}$}}
\def\eqn#1#2{\xdef #1{(\secsym\the\meqno)}\writedef{#1\leftbracket#1}%
\global\advance\meqno by1$$#2\eqno#1\eqlabeL#1$$}
%
\newskip\footskip\footskip14pt plus 1pt minus 1pt 
\def\footnotefont{\ninepoint}\def\f@t#1{\footnotefont #1\@foot}
\def\f@@t{\baselineskip\footskip\bgroup\footnotefont\aftergroup\@foot\let\next}
\setbox\strutbox=\hbox{\vrule height9.5pt depth4.5pt width0pt}
\global\newcount\ftno \global\ftno=0
\def\foot{\global\advance\ftno by1\footnote{$^{\the\ftno}$}}
%
\newwrite\ftfile
\def\footend{\def\foot{\global\advance\ftno by1\chardef\wfile=\ftfile
$^{\the\ftno}$\ifnum\ftno=1\immediate\openout\ftfile=foots.tmp\fi%
\immediate\write\ftfile{\noexpand\smallskip%
\noexpand\item{f\the\ftno:\ }\pctsign}\findarg}%
\def\footatend{\vfill\eject\immediate\closeout\ftfile{\parindent=20pt
\centerline{\bf Footnotes}\nobreak\bigskip\input foots.tmp }}}
\def\footatend{}
%
%
\global\newcount\refno \global\refno=1
\newwrite\rfile
%
\def\ref{\nref}
\def\nref#1{\xdef#1{[\the\refno]}\writedef{#1\leftbracket#1}%
\ifnum\refno=1\immediate\openout\rfile=refs.tmp\fi
\global\advance\refno by1\chardef\wfile=\rfile\immediate
\write\rfile{\noexpand\item{#1\ }\reflabeL{#1\hskip.31in}\pctsign}\findarg}
\def\findarg#1#{\begingroup\obeylines\newlinechar=`\^^M\pass@rg}
{\obeylines\gdef\pass@rg#1{\writ@line\relax #1^^M\hbox{}^^M}%
\gdef\writ@line#1^^M{\expandafter\toks0\expandafter{\striprel@x #1}%
\edef\next{\the\toks0}\ifx\next\em@rk\let\next=\endgroup\else\ifx\next\empty%
\else\immediate\write\wfile{\the\toks0}\fi\let\next=\writ@line\fi\next\relax}}
\def\striprel@x#1{} \def\em@rk{\hbox{}}
\def\lref{\begingroup\obeylines\lr@f}
\def\lr@f#1#2{\gdef#1{\ref#1{#2}}\endgroup\unskip}

\def\addref#1{\immediate\write\rfile{\noexpand\item{}#1}} 
\def\startrefs#1{\immediate\openout\rfile=refs.tmp\refno=#1}
\def\refs#1{\count255=1[\r@fs #1{\hbox{}}]}
\def\r@fs#1{\ifx\und@fined#1\message{reflabel \string#1 is undefined.}%
\nref#1{need to supply reference \string#1.}\fi%
\vphantom{\hphantom{#1}}\edef\next{#1}\ifx\next\em@rk\def\next{}%
\else\ifx\next#1\ifodd\count255\relax\xref#1\count255=0\fi%
\else#1\count255=1\fi\let\next=\r@fs\fi\next}
%

%
\newwrite\ffile\global\newcount\figno \global\figno=1
\def\fig{fig.~\the\figno\nfig}
\def\nfig#1{\xdef#1{fig.~\the\figno}%
\writedef{#1\leftbracket fig.\noexpand~\the\figno}%
\ifnum\figno=1\immediate\openout\ffile=figs.tmp\fi\chardef\wfile=\ffile%
\immediate\write\ffile{\noexpand\medskip\noexpand\item{Fig.\ \the\figno. }
\reflabeL{#1\hskip.55in}\pctsign}\global\advance\figno by1\findarg}
\def\xfig{\expandafter\xf@g}\def\xf@g fig.\penalty\@M\ {}
\def\figs#1{figs.~\f@gs #1{\hbox{}}}
\def\f@gs#1{\edef\next{#1}\ifx\next\em@rk\def\next{}\else
\ifx\next#1\xfig #1\else#1\fi\let\next=\f@gs\fi\next}
\newwrite\lfile
{\escapechar-1\xdef\pctsign{\string\%}\xdef\leftbracket{\string\{}
\xdef\rightbracket{\string\}}\xdef\numbersign{\string\#}}

\def\writestop{\def\writestoppt{\immediate\write\lfile{\string\pageno%
\the\pageno\string\startrefs\leftbracket\the\refno\rightbracket%
\string\def\string\secsym\leftbracket\secsym\rightbracket%
\string\secno\the\secno\string\meqno\the\meqno}\immediate\closeout\lfile}}
\def\writestoppt{}\def\writedef#1{}
\def\seclab#1{\xdef #1{\the\secno}\writedef{#1\leftbracket#1}\wrlabeL{#1=#1}}
\def\subseclab#1{\xdef #1{\secsym\the\subsecno}%
\writedef{#1\leftbracket#1}\wrlabeL{#1=#1}}
\newwrite\tfile \def\writetoca#1{}
\def\leaderfill{\leaders\hbox to 1em{\hss.\hss}\hfill}
\def\writetoc{\immediate\openout\tfile=toc.tmp
   \def\writetoca##1{{\edef\next{\write\tfile{\noindent ##1
   \string\leaderfill {\noexpand\number\pageno} \par}}\next}}}
%
%
%

\catcode`\@=12 
%
\edef\tfontsize{\ifx\answ\bigans scaled\magstep3\else scaled\magstep4\fi}
\font\titlerm=cmr10 \tfontsize \font\titlerms=cmr7 \tfontsize
\font\titlermss=cmr5 \tfontsize \font\titlei=cmmi10 \tfontsize
\font\titleis=cmmi7 \tfontsize \font\titleiss=cmmi5 \tfontsize
\font\titlesy=cmsy10 \tfontsize \font\titlesys=cmsy7 \tfontsize
\font\titlesyss=cmsy5 \tfontsize \font\titleit=cmti10 \tfontsize
\skewchar\titlei='177 \skewchar\titleis='177 \skewchar\titleiss='177
\skewchar\titlesy='60 \skewchar\titlesys='60 \skewchar\titlesyss='60
\def\titlefont{\def\rm{\fam0\titlerm}
\textfont0=\titlerm \scriptfont0=\titlerms \scriptscriptfont0=\titlermss
\textfont1=\titlei \scriptfont1=\titleis \scriptscriptfont1=\titleiss
\textfont2=\titlesy \scriptfont2=\titlesys \scriptscriptfont2=\titlesyss
\textfont\itfam=\titleit \def\it{\fam\itfam\titleit}\rm}
 \ifx\answ\bigans\else scaled\magstep1\fi
\ifx\answ\bigans\else

 \font\absi=cmmi10 scaled\magstep1
\font\absis=cmmi7 scaled\magstep1 \font\absiss=cmmi5 scaled\magstep1
\font\abssy=cmsy10 scaled\magstep1 \font\abssys=cmsy7 scaled\magstep1
\font\abssyss=cmsy5 scaled\magstep1 
\skewchar\absi='177 \skewchar\absis='177 \skewchar\absiss='177
\skewchar\abssy='60 \skewchar\abssys='60 \skewchar\abssyss='60
\fi
\font\ninerm=cmr9 \font\sixrm=cmr6 \font\ninei=cmmi9 \font\sixi=cmmi6
\font\ninesy=cmsy9 \font\sixsy=cmsy6 \font\ninebf=cmbx9
\font\nineit=cmti9 \font\ninesl=cmsl9 \skewchar\ninei='177
\skewchar\sixi='177 \skewchar\ninesy='60 \skewchar\sixsy='60
\def\ninepoint{\def\rm{\fam0\ninerm}
\textfont0=\ninerm \scriptfont0=\sixrm \scriptscriptfont0=\fiverm
\textfont1=\ninei \scriptfont1=\sixi \scriptscriptfont1=\fivei
\textfont2=\ninesy \scriptfont2=\sixsy \scriptscriptfont2=\fivesy
\textfont\itfam=\ninei \def\it{\fam\itfam\nineit}\def\sl{\fam\slfam\ninesl}%
\textfont\bffam=\ninebf \def\bf{\fam\bffam\ninebf}\rm}
%
%

\hyphenation{anom-aly anom-alies coun-ter-term coun-ter-terms}
\def\inv{^{\raise.15ex\hbox{${\scriptscriptstyle -}$}\kern-.05em 1}}

\def\Dsl{\,\raise.15ex\hbox{/}\mkern-13.5mu D} 
\def\dsl{\raise.15ex\hbox{/}\kern-.57em\partial}

\def\lspace{\ifx\answ\bigans{}\else\qquad\fi}
\def\lbspace{\ifx\answ\bigans{}\else\hskip-.2in\fi} 
\def\boxeqn#1{\vcenter{\vbox{\hrule\hbox{\vrule\kern3pt\vbox{\kern3pt
    \hbox{${\displaystyle #1}$}\kern3pt}\kern3pt\vrule}\hrule}}}
\def\mbox#1#2{\vcenter{\hrule \hbox{\vrule height#2in
        \kern#1in \vrule} \hrule}}  
%

\def\darr#1{\raise1.5ex\hbox{$\leftrightarrow$}\mkern-16.5mu #1}

\def\half{{\textstyle{1\over2}}} 
\def\roughly#1{\raise.3ex\hbox{$#1$\kern-.75em\lower1ex\hbox{$\sim$}}}

%
%


\def\frac#1#2{{#1\over#2}}

\def\half{\frac12}

\def\journal#1&#2(#3){\unskip, #1~\bf #2 \rm(19#3) }
\def\andjournal#1&#2(#3){\sl #1~\bf #2 \rm (19#3) }

\def\bra#1{\left\langle #1\right|}
\def\ket#1{\left| #1\right\rangle}

\def\One{{1\hskip -3pt {\rm l}}}
\catcode`\@=11\def\slash#1{\mathord{\mathpalette\c@ncel{#1}}}
\overfullrule=0pt
\def\steepslash{\c@ncel}
\def\frac#1#2{{#1\over #2}}

\def\:{\!:\!}
\def\inbar{\,\vrule height1.5ex width.4pt depth0pt}
\def\IQ{\relax\,\hbox{$\inbar\kern-.3em{\rm Q}$}}
\def\IB{\relax{\rm I\kern-.18em B}}
\def\IC{\relax\hbox{$\inbar\kern-.3em{\rm C}$}}
\def\IP{\relax{\rm I\kern-.18em P}}
\def\IR{\relax{\rm I\kern-.18em R}}
\def\ZZ{\relax\ifmmode\mathchoice
{\hbox{Z\kern-.4em Z}}{\hbox{Z\kern-.4em Z}}
{\lower.9pt\hbox{Z\kern-.4em Z}}
{\lower1.2pt\hbox{Z\kern-.4em Z}}\else{Z\kern-.4em Z}\fi}

\catcode`\@=12

\def\npb#1(#2)#3{{ Nucl. Phys. }{B#1} (#2) #3}
\def\plb#1(#2)#3{{ Phys. Lett. }{#1B} (#2) #3}
\def\pla#1(#2)#3{{ Phys. Lett. }{#1A} (#2) #3}
\def\prl#1(#2)#3{{ Phys. Rev. Lett. }{#1} (#2) #3}
\def\mpla#1(#2)#3{{ Mod. Phys. Lett. }{A#1} (#2) #3}
\def\ijmpa#1(#2)#3{{ Int. J. Mod. Phys. }{A#1} (#2) #3}
\def\cmp#1(#2)#3{{ Comm. Math. Phys. }{#1} (#2) #3}
\def\cqg#1(#2)#3{{ Class. Quantum Grav. }{#1} (#2) #3}
\def\jmp#1(#2)#3{{ J. Math. Phys. }{#1} (#2) #3}
\def\anp#1(#2)#3{{ Ann. Phys. }{#1} (#2) #3}
\def\prd#1(#2)#3{{ Phys. Rev. } {D{#1}} (#2) #3}
\def\ptp#1(#2)#3{{ Progr. Theor. Phys. }{#1} (#2) #3}
\def\aom#1(#2)#3{{ Ann. Math. }{#1} (#2) #3}

\def\bs{\bigskip}

\def\br{\buildrel}
\def\bra{\langle}
\def\ket{\rangle}

\def\B{{\bf B}}
\def\C{{\bf C}}

\def\P{{\bf P}}
\def\Q{{\bf Q}}
\def\R{{\bf R}}

\def\Z{{\bf Z}}

\def\cH{{\cal H}}

\def\cN{{\cal N}}

\def\cU{{\cal U}}
\def\cV{{\cal V}}

\input amssym
\def\ga{{\goth a}}
\def\gb{{\goth b}}
\def\gc{{\goth c}}
\def\gd{{\goth d}}

\def\gg{{\goth g}}
\def\gh{{\goth h}}

\def\gk{{\goth k}}
\def\gl{{\goth l}}

\def\gn{{\goth n}}

\def\gq{{\goth q}}

\def\gt{{\goth t}}
\def\gu{{\goth u}}

\def\da#1{{\p \over \p s_{#1}}}

\def\cicy#1(#2|#3)#4{\left(\matrix{#2}\right|\!\!
                     \left|\matrix{#3}\right)^{{#4}}_{#1}}

\def\Lra{\Longrightarrow}
\def\LRa{\Leftrightarrow}
\def\hra{\hookrightarrow}
\def\hla{\hookleftarrow}
\def\ra{\rightarrow}
\def\la{\leftarrow}
\def\thra{\twoheadrightarrow}

\def\da{\downarrow}
\def\ua{\uparrow}

\def\bs{\bigskip}

\def\Box{{\,\lower0.9pt\vbox{\hrule
\hbox{\vrule height 0.2 cm \hskip 0.2 cm
\vrule height 0.2 cm}\hrule}\,}}

\global\newcount\thmno \global\thmno=0
\def\definition#1{\global\advance\thmno by1
\bigskip\noindent{\bf Definition \secsym\the\thmno. }{\it #1}
\par\nobreak\medskip\nobreak}
\def\question#1{\global\advance\thmno by1
\bigskip\noindent{\bf Question \secsym\the\thmno. }{\it #1}
\par\nobreak\medskip\nobreak}
\def\theorem#1{\global\advance\thmno by1
\bigskip\noindent{\bf Theorem \secsym\the\thmno. }{\it #1}
\par\nobreak\medskip\nobreak}
\def\proposition#1{\global\advance\thmno by1
\bigskip\noindent{\bf Proposition \secsym\the\thmno. }{\it #1}
\par\nobreak\medskip\nobreak}
\def\corollary#1{\global\advance\thmno by1
\bigskip\noindent{\bf Corollary \secsym\the\thmno. }{\it #1}
\par\nobreak\medskip\nobreak}
\def\lemma#1{\global\advance\thmno by1
\bigskip\noindent{\bf Lemma \secsym\the\thmno. }{\it #1}
\par\nobreak\medskip\nobreak}
\def\conjecture#1{\global\advance\thmno by1
\bigskip\noindent{\bf Conjecture \secsym\the\thmno. }{\it #1}
\par\nobreak\medskip\nobreak}
\def\exercise#1{\global\advance\thmno by1
\bigskip\noindent{\bf Exercise \secsym\the\thmno. }{\it #1}
\par\nobreak\medskip\nobreak}
\def\remark#1{\global\advance\thmno by1
\bigskip\noindent{\bf Remark \secsym\the\thmno. }{\it #1}
\par\nobreak\medskip\nobreak}
\def\problem#1{\global\advance\thmno by1
\bigskip\noindent{\bf Problem \secsym\the\thmno. }{\it #1}
\par\nobreak\medskip\nobreak}
\def\others#1#2{\global\advance\thmno by1
\bigskip\noindent{\bf #1 \secsym\the\thmno. }{\it #2}
\par\nobreak\medskip\nobreak}
\def\proof{\noindent Proof: }

\def\thmlab#1{\xdef #1{\secsym\the\thmno}\writedef{#1\leftbracket#1}\wrlabeL{#1=#1}}
%
%
\def\newsec#1{\global\advance\secno by1\message{(\the\secno. #1)}
\global\subsecno=0\thmno=0\eqnres@t\noindent{\bf\the\secno. #1}
\writetoca{{\secsym} {#1}}\par\nobreak\medskip\nobreak}
\def\eqnres@t{\xdef\secsym{\the\secno.}\global\meqno=1\bigbreak\bigskip}
\def\sequentialequations{\def\eqnres@t{\bigbreak}}\xdef\secsym{}
%

%
\newcount{\exnum}
\def\prob{\advance\exnum by 1
\bigskip\item{\bf\the\exnum.}\ }
\newcount{\exnum}
\def\next{\advance\exnum by 1
\bigskip\noindent{\the\exnum.}\ }
\def\np{\vfill\eject}